\documentclass[12pt]{article} %

\usepackage{graphicx}
\usepackage{epstopdf}
\usepackage{amssymb}
\usepackage{amsmath}
\usepackage{dsfont}
\usepackage{bm}

\usepackage{lscape}
\usepackage{color}

\newcommand{\eod}{{\hfill $\square$}}

\usepackage{picins}

\usepackage{multirow}

\newtheorem{theorem}{Theorem}
\newtheorem{coro}{Corollary}

\newtheorem{remark}{Remark}

\newtheorem{defi}{Definition}

\usepackage{balance}
\makeatletter
\def\blfootnote{\xdef\@thefnmark{}\@footnotetext}
\makeatother

\usepackage{keyval}
\makeatletter
\define@key{cascading}{width}{\cascading@wd=#1}
\define@key{cascading}{sep}{\def\cascading@sep{#1}}
\newdimen\cascading@wd
\newcommand{\cascadingblocks}[2][]{
  \setkeys{cascading}{sep=2ex,#1}%
  \leavevmode\vbox{\offinterlineskip
    \@for\next:=#2\do{%
      \settowidth{\dimen@}{\next}%
      \ifdim\dimen@>\cascading@wd
        \cascading@wd=\dimen@
      \fi
    }%
    \@for\next:=#2\do{%
      \cascading@rule
      \hbox{\fbox{\hb@xt@\cascading@wd{\hss\next\hss}}}%
    }
  }
}
\def\cascading@rule{%
  \def\cascading@rule{%
    \hb@xt@\dimexpr\cascading@wd+2\fboxsep+2\fboxrule\relax
      {\hss\vrule\@height\cascading@sep\hss}%
  }%
}
\makeatother

\usepackage{smartdiagram}
\begin{document}
\title{Matrix-Valued  Mean-Field-Type Games:  Risk-Sensitive, Adversarial, and Risk-Neutral Linear-Quadratic Case}
\date{First draft: September 2018. This version: June  2019}
\author{Julian Barreiro-Gomez, Tyrone E. Duncan  and Hamidou Tembine
\thanks{Julian Barreiro-Gomez and Hamidou Tembine are with Learning \& Game Theory Laboratory, New York University Abu Dhabi (NYUAD), Saadiyat Campus PO Box 129188, United Arab Emirates, (e-mails: jbarreiro@nyu.edu, tembine@nyu.edu).}
\thanks{Tyrone Duncan is with Department of Mathematics, University of Kansas, Lawrence, KS 66044, USA, (e-mail: duncan@ku.edu).}
\thanks{We gratefully acknowledge support from U.S. Air Force Office of Scientific Research under grant number FA9550-17-1-0259 and FA9550-17-1-0073 and National Science Foundation under grant number DMS 1411412.}
}

\maketitle

\newpage

\begin{abstract}
In this paper we study a class of matrix-valued linear-quadratic mean-field-type games for both the risk-neutral, risk-sensitive and robust cases. Non-cooperation, full cooperation and adversarial between teams are treated.  We provide a semi-explicit solution for both problems by means of a direct method. 
The state dynamics is described by a matrix-valued linear jump- diffusion-regime switching system of  conditional mean-field type where the conditioning is with respect to common noise which is a regime switching process.  The optimal strategies are in a state-and-conditional mean-field feedback form. Semi-explicit solutions of equilibrium costs and strategies  are also provided for the full cooperative, adversarial teams, risk-sensitive full cooperative and risk-sensitive adversarial team cases. It is shown that   full cooperation increases the  well-posedness  domain under risk-sensitive  decision-makers by means of population risk-sharing. Finally, relationship between risk-sensitivity and  robustness are established in the mean-field type context.
\end{abstract}

\tableofcontents

\section{Introduction}

The Markowitz paradigm, also termed as  the mean-variance paradigm, is often characterized as dealing with portfolio risk and (expected) return \cite{markowitz1,markowitz2,markowitz3}. A typical example of risk concerns in the current online market is the evolution of   prices of the digital and cryptocurrencies (bitcoin, litecoin, ethereum, dash, and other altcoins etc).  Variance plays a base model for many risk measures.  Here, we address  variance reduction problems when several decision-making entities are involved. When the decisions made by the entities influence each other, the decision-making is said to be interactive (interdependent). Such problems are termed as game problems.  Game problems in which the state dynamics is given  a linear stochastic system with a Brownian motion and a cost functional that is quadratic in the state and the control,  is often called the linear-quadratic Gaussian (LQG) games. 
For  generic LQG game problems under perfect state observation,  the optimal strategy of the decision-maker is a linear state-feedback strategy which is identical to an optimal control for the corresponding deterministic linear-quadratic game problem where the Brownian motion is replaced by the zero process. Moreover the equilibrium cost only differs from the deterministic game 
problem's equilibrium cost by the integral of a function of time. However, when the diffusion (volatility) coefficient is state and control-dependent, the structure of the resulting differential system  as well as the equilibrium cost vector are modified.   These results were widely known in both dynamic optimization, control and game theory literature. 
For both  LQG control and LQG zero-sum  games, it can be shown that a simple square completion method, provides an explicit solution to the problem. It was successfully  developed and applied by 
Duncan et al. \cite{new1,new2,ref0,ref1,ref2,ref4} in the mean-field-free case. Moreover,  Duncan et al.  have extended the direct method to more general noises including fractional Brownian noises and some non-quadratic cost functionals on spheres and torus. 
Inspired by applications in engineering (internet connexion, battery state etc) and in finance (price, stock option, multi-currency exchange etc) were not only Gaussians but also jump processes (Poisson, L\'evy etc) play important features, the question of  extending the framework to linear-quadratic games under state dynamics driven by jump-diffusion processes were naturally posed.  Adding a  Poisson jump and regime switching may allow to capture in particular larger jumps which may not be captured by just increasing diffusion coefficients. Several examples such as multi-currency exchange or  cloud-server rate allocation on blockchains are naturally in matrix form.

The main goal of this work is to investigate whether Direct Method can be used to solve  matrix-valued risk-sensitive and adversarial robust  mean-field-type game problems which are non-standard problems \cite{alain2}. To do so,
we  modify the state dynamics to include  conditional mean-field terms which are 
\begin{itemize}
\item the conditional expectation of the matrix-valued state with respect to the filtration of the common noise which is a regime switching process, is added to   the drift, diffusion, and jump coefficient functionals.
\item the  conditional expectation of the matrix-valued control-actions, is included in the drift, diffusion, jump coefficient functional.
\end{itemize}
 We also modify the instant cost and terminal cost function to include  
 \begin{itemize} \item the square of the conditional expectation of the matrix-valued state and 
 \item  the square of the conditional expectation of the matrix-valued control action.  
 \end{itemize} Involving these features lead to matrix-valued
mean-field-type game theory which focuses on 
{\it (matrix-valued) games with distribution-dependent quantity-of-interest} such as payoff, cost and state dynamics. 
It can be seen as the multiple agent generalization of single agent (matrix-valued) mean-field-type control problem \cite{alain2}.

If the state dynamics and or the cost functional involve a mean-field term (such as the  expectation  of the state and or   the expectation of the control actions), the game is said to be a  LQG game of mean-field type, or MFT-LQG.
 For such game problems, various solution methods such as the Stochastic Maximum Principle (SMP) (\cite{alain2}) and the  Dynamic Programming Principle (DPP) with Hamilton-Jacobi-Bellman-Isaacs equation and Fokker-Planck-Kolmogorov equation have been proposed \cite{ref00,alain2,alaintem}. 
 
If the state dynamics and or the cost functional involve a conditional mean-field term (such as the conditional expectation  of the matrix-valued state and or   the conditional expectation of the matrix-valued control actions), the game is said to be a  matrix-valued LQG game of conditional mean-field type, or cMFT-LQG (or conditional McKean-Vlasov matrix-valued LQG games).  If in addition, the matrix-valued state dynamics is driven by a matrix-valued jump-diffusion process, then the problem is termed as a cMFT-LQJD matrix-valued game problem.  
We aim to study the behavior of such cMFT-LQJD matrix-valued game problems  when conditional mean-field terms are involved.

Games with global uncertainty and common noise have been widely suggested in the literature.
Anonymous sequential and mean-field games with common noise can be considered as a natural generalization of the mean-field game problems  (see \cite{jova,jova2} and the references therein). 
The works in \cite{bergin1,bergin2}  considered mean-field games with common noise and 
obtained optimality system that determine mean-field equilibria conditioned of the information.
The work in \cite{lions18,ref6} provides sufficiency conditions for well-posedness of mean-field games with common noise and a major player. Existence of solutions of the resulting stochastic optimality systems are examined in \cite{ref4ty}. A probabilistic approach to the master equation is developed in \cite{ref5}.  In order to determine the optimal strategies of the decision-maker, the previous works used a maximum principle or a master equation which involves a stochastic Fokker-Planck equation (see \cite{ref6,ref8ty,ref9,ref7} and the references therein).

 Most studies illustrated mean-field game methods in the linear-quadratic game with infinite number of decision-makers \cite{bardi1,bardi2,Boualem2014v2,robust,TAC2014}.  These works assume indistinguishability within classes and the cost functionals were assumed to be identical or invariant per permutation of decision-makers indexes. Note that the indistinguishability assumption is not fulfilled for many interesting problems such as  variance reduction or 
and  risk quantification problems in which decision-makers have different sensitivity towards the risk. One typical and practical example is to consider a multi-level building in which every resident has its own comfort zone temperature and aims to use the Heating, Ventilating, and Air Conditioning  (HVAC) system to be in  its comfort temperature zone and maintain it within its own comfort zone. This problem clearly does not satisfy the indistinguishability assumption used in the previous works on mean-field games. Therefore, it is reasonable to look at the problem beyond the indistinguishability assumption.
 Here we drop these assumptions and dealt with the problem directly with arbitrarily  finite number of decision-makers.
 In the LQ-mean-field game problems the state process can be modeled by a set of linear 
stochastic differential equations of McKean-Vlasov and the preferences are formalized by quadratic or exponential of integral of quadratic cost functions with mean-field terms. These game  problems are of practical interests  and a
detailed exposition of this theory can be found in \cite{alain2,alain,ref3,ref3ty,Engwerd,ref8}. The popularity of these game problems is due to practical considerations
 in consensus problems, signal processing, pattern recognition, filtering, prediction, economics and management science \cite{ref8ty,Boualem2014,book2,rt4}.  

To some extent, most of the risk-neutral versions of these optimal controls are analytically and numerically solvable \cite{ref1,ref3,ref4,new1,new2}. On the other hand, the linear quadratic robust  setting naturally appears if the decision makers' objective is to minimize the effect of a small perturbation and related variance of the optimally controlled nonlinear process. By solving a linear quadratic   game problem of mean-field type, and using the implied optimal control actions, decision-makers can significantly reduce the variance (and the cost) incurred by this perturbation.  The variance reduction and minmax problems have very interesting applications in risk quantification problems under adversarial attacks and in security issues in interdependent infrastructures and networks \cite{rt2,rt1,rt3,rt4,wireless}. In \cite{Djehiche_2017}, the control for the evacuation of a multi-level building is designed by means of mean-field games and mean-field-type control.
In \cite{Djehiche_2018}, electricity price dynamics in the smart grid is analyzed  using a mean-field-type game approach under common noise which is of diffusion type. Risk-neutral Linear-Quadratic MFT-LQJD games have been studied for the one dimensional case in \cite{BaDuTe_2018}.

\subsection*{Our contribution}
In this paper, we use a simple argument that gives the  risk-neutral equilibrium strategy and  robust adversarial  mean-field-type saddle point for a class of cMFT-LQJD  matrix-valued games without use of the well-known solution 
methods (SMP and  DPP). We apply a basic It\^o's formula following by a square completion method in the  risk-neutral/adversarial mean-field-type matrix-valued game problems. It is shown that this method is well suited to cMFT-LQJD risk-neutral/robust games as well as to variance  reduction  performance functionals with jump-diffusion-regime switching common noise. Applying the solution methodology related to the DPP or the SMP requires involved (stochastic) analysis and convexity arguments to generate necessary and sufficient optimality criteria. We avoid all this with this method. 

Zero-sum stochastic differential games are important class of stochastic games. The optimality system leads to an
Hamilton-Jacobi-Bellman-Isaacs (HJBI) system of equations, which is an extension of the HJB equation to stochastic differential games. When common noise is involved it becomes a stochastic HJBI system. Obviously, studying well-posedness, existence and uniqueness of such is a challenging task because of the minmax  and  maxmin operators. Usually, upper value and lower value equilibrium payoffs are investigated. In addition, when conditional mean-field terms are involved as it is the case here, the system is coupled with a stochastic Fokker-Planck-Kolmogorov system leading to a master system. Here we provide an easy way to solve such a system of means of a direct method. 

Relationship between risk-sensitive and roust  conditional mean-field-type games are established in the case without jump and with a single regime.

Our contribution can be summarized as follows. 
We formulate and solve a matrix-valued linear-quadratic mean-field-type game described by a linear jump-diffusion dynamics and a mean-field-dependent quadratic cost functional that is conditioned a common noise which includes not only a Brownian motion but also a jump process and regime switching. 
Since the matrices are switching dependent, they can be seen as random coefficients. 
The optimal strategies for the decision-makers are given semi-explicitly using a simple and direct method based on square completion, suggested in Duncan {\it et al.} in e.g.  \cite{ref0} for the mean-field free case. This approach does not use the well-known solution methods such as the Stochastic Maximum Principle and the Dynamic Programming Principle with  stochastic Hamilton-Jacobi-Bellman-Isaacs equation and stochastic Fokker-Planck-Kolmogorov equation. It does require  extended stochastic backward-forward partial  integro-differential equations (PIDE) to solve the problem. 
 In the risk-neutral   linear-quadratic mean-field-type game with perfect state observation and with common noise, we show that, generically  there is a minmax strategy to the conditional mean of the state and provide a  sufficient condition of existence  of mean-field-type saddle point.   Sufficient conditions for existence and uniqueness of robust mean-field equilibria are obtained when the horizon length   is small enough and the Riccati coefficients are almost surely positive.

\begin{figure}[t!]
\centering \includegraphics[width=\columnwidth]{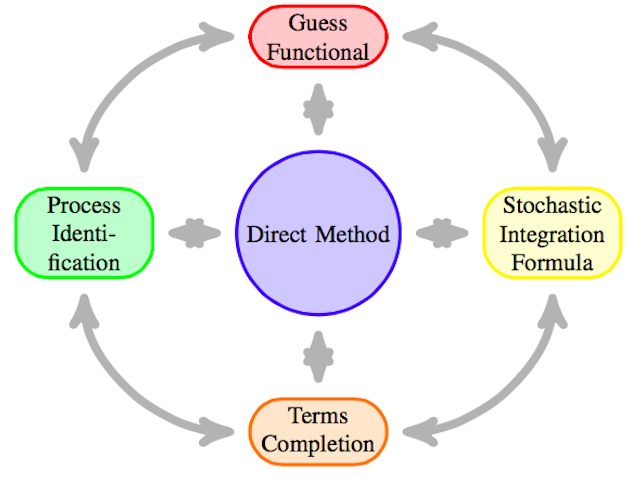}\\
\label{fig:sum:directad}  \caption{Direct method and its key   steps.}
  \end{figure}

 In addition, this work extends the results in \cite{Duncannew} in various ways: 
 \begin{itemize}
 \item Extension to matrix-form of arbitrary dimensions.
 \item the common noise which is  a regime switching was not considered in \cite{Duncannew}. 

 \item The solution here involves a matrix-valued  differential system which differs from the results in  \cite{Duncannew}.
 
  \end{itemize}
To solve the aforementioned problem in a semi-explicit way, we follow a direct method. The method starts by identifying a partial guess functional where the coefficient functionals are random and regime switching dependent. Then, it uses It\^o's formula for jump-diffusion-regime switching processes, followed by a completion of squares for both control  and conditional mean control. Finally, the processes are identified using an orthogonal decomposition technique and stochastic  differential equations are derived in a semi-explicit way. The procedure is summarized in Figure \ref{fig:sum:directad}. The contributions of this work are summarized in Table \ref{tab:resultsad}. 

 \begin{table}[ht!]
	\hspace{-4cm}
  		\begin{center}
  		\begin{tabular}{ccc}
  			\hline
  			& \multicolumn{2}{c}{\textbf{Reference}}\\
			
  			\hline
  			\textbf{Cost function Terms} & \cite{Duncannew}    & this work  \\
  			\hline
			
			risk-neural non-cooperation&   $\checkmark$, 1D  &  $\checkmark$, matrix-valued \\
			risk-neural full-cooperation &   $\checkmark$, 1D  &  $\checkmark$, matrix-valued \\
			risk-neural adversarial/robust &  $\checkmark$, 1D  &  $\checkmark$, matrix-valued \\
			risk-sensitive non-cooperation&     &  $\checkmark$, matrix-valued \\
			risk-sensitive full-cooperation &    &  $\checkmark$, matrix-valued \\
			risk-sensitive adversarial/robust &     &  $\checkmark$, matrix-valued \\

			\hline
			\textbf{Drift Terms} &   &   \\
  			\hline

  			$B_1,B_2$ &   $\checkmark$  &  $\checkmark$ \\
  			$\bar{B}_1,\bar{B}_2$ &  $\checkmark$ &  $\checkmark$ \\
  			
  			\hline
  			\textbf{Diffusion Terms} &  &    \\
  			\hline
  			$S_0$ &  $\checkmark$  & $\checkmark$ \\

  			\hline
  			\textbf{Jump Terms} &  &      \\
  			\hline
  			$M_{0}$ &  $\checkmark$&  $\checkmark$  \\
  		
  			\hline
  			\textbf{Switching Term} &  &     \\
  			\hline
  			$s$ &  &    $\checkmark$ \\
			
  			\hline
  			\textbf{Cost function Terms} &  &      \\
  			\hline
  			$\mathbb{E}[X],\mathbb{E}[U]$ &  $\checkmark$ &     $\checkmark$ \\
  			$\mathrm{var}[X],\mathrm{var}[U]$ & $\checkmark$ & $\checkmark$ \\
  			  			\hline
%
  			\textbf{Number of Decision makers} &  &      \\
  			\hline
  			One single team & $\checkmark$  & $\checkmark$  \\
  			Multiple players & $\checkmark$ &   $\checkmark$  \\
			Two adversaries & $\checkmark$ &   $\checkmark$  \\
			Two adversarial teams&  &   $\checkmark$  \\
  			\hline
%
  		\end{tabular}
  	\end{center} \caption{Contributions  with respect to recent literature.}
  	\label{tab:resultsad}
  \end{table}

{\it To  the best of the authors knowledge this is the first work to consider regime switching in matrix-valued mean-field-type game theory.}
%
  
 A brief outline of the rest of the paper follows. The next section introduces a generic  game model. After that,  the MF-LQJD conditional mean-field-type game  problem is investigated and its solution is presented.  The last section concludes the paper.  

 \subsection*{Notation and Preliminaries}

We introduce the following notations. 
Let $[0,T], \ T>0$  be a fixed time horizon and $(\Omega,\mathcal{F},\mathbb{F}^{B,N,s}, \mathbb{P})$  be a given filtered
probability space.The filtration $\mathbb{F}=\{{\mathcal{F}}^{B,N,s}_t,\ 0\leq t \leq T\}$ is the natural filtration of the union $\{B,N,s\}$ augmented by $\mathbb{P}-$null sets of ${\mathcal F}.$ 
In practice, $B$ is used to capture smaller disturbance and $N$ is used for larger jumps of the system. 


An admissible control strategy $U_i$ of the decision-maker $i$  is an $\mathbb{F}$-adapted and square-integrable process with values in a 
 $\mathbb{U}_i=\mathbb{R}^{d\times d}, \ d\geq 1$. We denote the set of all admissible controls
 by $\mathcal{U}_i$.

\section{Problem Formulation}
We consider $n\geq 1$ decision-makers over the time  horizon $[0,T], \ T>0.$   Each decision-maker $i$ chooses a matrix-valued strategy $U_i$ over the horizon $[0,T].$ The state satisfies the following matrix-valued linear jump-diffusion-regime switching system of mean-field type:

\begin{equation} \label{equation:state}
\begin{array}{l}
dX = [ B_1(X - \bar{X}) + (B_1 + \bar{B}_1)\bar{X} \\
+ \sum_{j=1}^{n} B_{2j} (U_j - \bar{U}_j) + \sum_{j=1}^{n} (B_{2j}+\bar{B}_{2j}) \bar{U}_j ]dt\\
+ S_0 dB + \int_{\theta} M_0(\cdot,\theta) \tilde{N}(dt,d\theta),
\end{array}
\end{equation}
where 
$X(t)$, $\bar{X}(t)= \mathbb{E}[X(t)|\mathcal{F}^s_t]$,  $B_1(t,s(t))$, $\bar{B}_1(t,s(t)) \in \mathbb{R}^{d\times d},$
$ B_{2j}(t,s(t))$, $\bar{B}_{2j}(t,s(t))$,  $S_0(t,s(t))$, $M_0(t,s(t)) \in \mathbb{R}^{d\times d},$ $B$ is a $\mathbb{R}^{d\times d}$ matrix-valued Brownian motion,  $U_i \in \mathcal{U}_i$, 
$s$ is a regime switching process with transition rates $\tilde{q}_{ss'}$ satisfying $\tilde{q}_{ss'}>0$, $\tilde{q}_{ss}=-\sum_{s'\ne s}\tilde{q}_{ss'}.$.
$N(t,.)$ is a $\mathbb{R}^{d\times d}$ matrix-valued Poisson random process with compensated process  $\tilde{N}(dt,d\theta) =N(dt,d\theta)- \nu(d\theta)dt,$  $\nu$ is a matrix of Radon measure over the set of jump sizes $\Theta.$  $\mathcal{F}^s_t$ is the filtration generated by the regime switching process $s.$ By abuse of notation we omit the  use of $s(t-)$ and $X(t-)$ for the left value of switching $s(t)$ and  the jump $X(t)$ respectively. 

To the state system (\ref{equation:state}), we associate the cost functional of decision-maker $i$ 

\begin{align}
\label{eq:main_cost}
\begin{array}{l}
L_i(X,s,U)  
= \langle Q_i(T,s(T)) (X(T) - \bar{X}(T)), X(T) - \bar{X}(T)  \rangle \\
+ \langle (Q_i(T,s(T))+\bar{Q}_i(T,s(T))) \bar{X}(T), \bar{X}(T)  \rangle\\
+ \int_{0}^{T} \langle Q_i (X - \bar{X}), X - \bar{X}  \rangle + \langle (Q_i+\bar{Q}_i) \bar{X}, \bar{X}  \rangle dt\\
+ \int_{0}^{T} \langle R_i (U_i - \bar{U}_i), U_i - \bar{U}_i  \rangle + \langle (R_i+\bar{R}_i) \bar{U}_i, \bar{U}_i  \rangle dt,
\end{array}
\end{align}
where $\langle A, B\rangle=\mbox{trace}(A^*B)=\mbox{trace}(B^*A), $  $A^*$  being the adjoint operator of $A$ (transposition),  The coefficients  $Q_i:=Q_i(t,s(t)),\ R_i:=R_i(t,s(t)), \bar{Q}_i:=\bar{Q}_i(t,s(t)), \bar{R}_i:=\bar{R}_i(t,s(t))$ are possibly time and regime-switching dependent with values in  $\mathbb{R}^{d\times d}.$

{\color{blue} 
The reader may want to know why a matrix-valued dynamics instead of a vector-valued dynamics. One typical example is the evolution of the rate of change of blockchain tokens and classical currencies. The exchange rate between token $t_k$ and token $t_{k'}$ is given by the entries $(k,k')$ of $X$ with value $X_{kk'}.$ Since these tokens are correlated one obtains a matrix-valued process $X.$ 
}

\subsection{Risk-Neutral}
We provide basic definitions of the risk-neutral problems and their solution concepts.  
\begin{defi} [Mean-Field-Type Risk-Neutral Best-Response]
Given $(U_j, \ j\neq i),$ a risk-neutral best response strategy of decision-maker $i$ is a strategy that solves
$ \inf_{U_i} \ \mathbb{E}L_i \ $ subject to  (\ref{equation:state}). The set of risk-neutral best responses of $i$ is denoted by $\mbox{rnBR}_i (U_{-i}).$ \eod
\end{defi}

\begin{defi} [Mean-Field-Type Risk-Neutral Nash Equilibrium]
A mean-field-type risk-neutral Nash equilibrium  is a strategy profile $(U^\mathrm{rn}_j, \ j\in \{1,\ldots, n\}),$ of all decision-makers such that for every decision-maker $i,$  $$U_i^\mathrm{rn}\in \mbox{rnBR}_i (U^\mathrm{rn}_{-i}).$$ \eod
\end{defi}
\begin{defi} [Mean-Field-Type Risk-Neutral Full Cooperation] A mean-field-type risk-neutral fully cooperative solution  is a strategy profile $(U^\mathrm{rn,g}_j, \ j\in \{1,\ldots, n\}),$ of all decision-makers such that 
$$\mathbb{E}[L_0(X,s,U^{rn,g})]=\inf_{(U_1, \ldots, U_n)}\mathbb{E}[L_0(X,s,U],$$ where $$L_0:=\sum_{j=1}^n L_j,$$ is the social (global) cost. \eod
\end{defi}

\begin{defi} [Mean-Field-Type Risk-Neutral Saddle Point Solution ] The set of decision-makers is divided into two teams. A team of defenders and a team of attackers. The defenders set is
 $$I_{+}:=\{ i\in \{1,\ldots, n\} | \ R_i\succ 0,  (R_i+\bar{R}_i) \succ 0\}$$ and the attackers set is $$I_{-}:=\{ j\in \{1,\ldots, n\} | \ -R_j\succ 0, \  -(R_j+\bar{R}_j) \succ 0\}.$$ 
A mean-field-type risk-neutral saddle point   is a strategy profile $(U^\mathrm{ad}_j, \ j\in I_{+}),$ of the team of defenders and  $(U^\mathrm{ad}_j, \ j\in I_{-})$ of the team of attackers such that 
\begin{align*}
\begin{array}{l}
L^\mathrm{ad}(X,s,(U^\mathrm{ad}_i)_{i\in I_{+} }, (U_j)_{j\in I_{-} })  \leq L^\mathrm{ad}(X,s,U^\mathrm{ad}) \\
\leq  L^\mathrm{ad}(X,s,(U_i)_{i\in I_{+} }, (U^\mathrm{ad}_j)_{j\in I_{-} })
\end{array}
\end{align*}
and $L^\mathrm{ad}(X,s,U^\mathrm{ad})$ is the value of the adversarial team  (risk-neutral) game, where 
\begin{align*}
\begin{array}{l}
L^\mathrm{ad}(X,s,U) \\
:= \langle Q(T,s(T)) (X(T) - \bar{X}(T)), X(T) - \bar{X}(T)  \rangle \\
+ \langle (Q(T,s(T))+\bar{Q}(T,s(T))) \bar{X}(T), \bar{X}(T)  \rangle\\
+ \int_{0}^{T} \langle Q (X - \bar{X}), X - \bar{X}  \rangle + \langle (Q+\bar{Q}) \bar{X}, \bar{X}  \rangle dt\\
+ \int_{0}^{T} \sum_{i=1}^n\langle R_i (U_i - \bar{U}_i), U_i - \bar{U}_i  \rangle \\
+ \sum_{i=1}^n \langle (R_i+\bar{R}_i) \bar{U}_i, \bar{U}_i  \rangle dt.
\end{array}
\end{align*}
\eod
\end{defi}
\subsection{Risk-Sensitive}
We provide basic definitions of  risk-sensitive problems and their solution concepts.  
\begin{defi} [Mean-Field-Type Risk-Sensitive Best-Response]
Given $(U_j, \ j\neq i),$ a risk-sensitive best response strategy of decision-maker $i$ is a strategy that solves
$$ \inf_{U_i} \ \frac{1}{\lambda_i}\log \left(\mathbb{E} [e^{\lambda_i L_i}]\right), \ $$ subject to  (\ref{equation:state}). 
The set of risk-sensitive best responses of $i$ is denoted by $\mbox{rsBR}_i (U_{-i}).$ \eod
\end{defi} 

For $\lambda_i\neq 0,$ the risk-sensitive loss functional $$ \frac{1}{\lambda_i}\log \left(\mathbb{E} [e^{\lambda_i L_i}]\right)$$  includes not only the first moment $\mathbb{E}[L_i]$ but also all the higher moments $\mathbb{E}[L^k_i], \ k\geq 1.$ 

\begin{defi} [Mean-Field-Type Risk-Sensitive Nash Equilibrium ]
A mean-field-type risk-sensitive Nash equilibrium  is a strategy profile $(U^\mathrm{rs}_j, \ j\in \{1,\ldots, n\}),$ of all decision-makers such that for every decision-maker $i,$  $$U_i^\mathrm{rs}\in \mbox{rsBR}_i (U^\mathrm{rs}_{-i}).$$ \eod
\end{defi}

\begin{defi} [Mean-Field-Type Risk-Sensitive Full Cooperation] A mean-field-type risk-sensitive fully cooperative solution  is a strategy profile $(U^\mathrm{rs,g}_j, \ j\in \{1,\ldots, n\}),$ of all decision-makers such that 
$$\inf_{(U_1, \ldots, U_n)} \frac{1}{\lambda}\log \left[ \mathbb{E}e^{\lambda L_0(X,s,U)} \right]= \frac{1}{\lambda}\log \left[ \mathbb{E}e^{\lambda L_0(X,s,U^{rs,g})} \right],$$ 
\end{defi}

\begin{defi} [Mean-Field-Type Risk-Sensitive Saddle Point Solution ] The set of decision-makers is divided into two teams. A team of defenders and a team of attackers. The defenders set is $I_{+}:=\{ i\in \{1,\ldots, n\} | \ R_i\succ 0,  (R_i+\bar{R}_i) \succ 0\}$ and the attackers set is $I_{-}:=\{ j\in \{1,\ldots, n\} | \ -R_j\succ 0, \  -(R_j+\bar{R}_j) \succ 0\}.$ 
A mean-field-type risk-sensitive saddle point   is a strategy profile $(U^\mathrm{ad}_j, \ j\in I_{+}),$ of the team of defenders and  $(U^\mathrm{ad}_j, \ j\in I_{-})$ of the team of attackers such that 
\begin{align*}
\begin{array}{l}
\frac{1}{\lambda}\log \left[ \mathbb{E}e^{\lambda L^\mathrm{ad}(X,s,(U^\mathrm{ad}_i)_{i\in I_{+} }, (U_j)_{j\in I_{-} })} \right]  \\
\leq \frac{1}{\lambda}\log \left[ \mathbb{E}e^{\lambda L^\mathrm{ad}(X,s,U^\mathrm{ad})}\right] \\
\leq  \frac{1}{\lambda}\log \left[ \mathbb{E}e^{\lambda L^\mathrm{ad}(X,s,(U_i)_{i\in I_{+} }, (U^\mathrm{ad}_j)_{j\in I_{-} })}\right].
\end{array}
\end{align*}
\end{defi}

\section{Main Results }  
This section presents the main results of the article.

\subsection{Risk-Neutral Case }
We start with the risk-neutral Nash equilibrium problem.
\begin{theorem} \label{thm:1}
Assume that  $Q_i, R_i, Q_i+\bar{Q}_i, R_i+\bar{R}_i$ are symmetric positive definite.  
Then the matrix-valued mean-field-type (risk-neutral) Nash equilibrium strategy and the (risk-neutral) equilibrium  cost are given by:
\begin{align*} 
\left\{
\begin{array}{l}
U^\mathrm{rn}_i - \bar{U}^\mathrm{rn}_i =- \frac{1}{2} R_i^{-1}B_{2i}^* (P^*_i+P_i)(X - \bar{X}),\\
\bar{U}^\mathrm{rn}_i =  -\frac{1}{2} (R_i+\bar{R}_i)^{-1} (B_{2i}+\bar{B}_{2i})^* (\bar{P}^*_i+\bar{P}_i)\bar{X},\\
L_i^\mathrm{rn}(X,s,U^\mathrm{rn}) = \mathbb{E}\langle P_i(0,s(0)) (X_0 - \bar{X}_0), X_0 - \bar{X}_0  \rangle \\
+ \mathbb{E}\langle \bar{P}_i(0,s(0))  \bar{X}_0, \bar{X}_0  \rangle  + \mathbb{E}[\delta_i(0,s(0))],\\
i \in \{1,\dots, n\},
\end{array}
\right.
\end{align*} 
where $P_i, \bar{P}_i,$ and $\delta_i$ solve the following differential equations:
\begin{equation} \label{tt0}
\left\{
\begin{array}{l}
\dot{P}_i + Q_i +P_i B_1+B_1^*P_i   + \sum_{s'\neq s}({P}_i(t,s')-{P}_i(t,s))\tilde{q}_{ss'} \\
- \frac{1}{4} (P^*_i+P_i) B_{2i} R_i^{-\frac{1}{2}*} R_i^{-\frac{1}{2}} B_{2i}^* (P^*_i+P_i) \\
{\color{blue} - \frac{1}{4} \sum_{j \neq i} (P^*_j+P_j)B_{2j} R_j^{-1*}B_{2j}^*(P^*_i+P_i) }\\
{\color{blue} - \frac{1}{4} \sum_{j \neq i} (P^*_i+P_i)B_{2j} R_j^{-1*}B_{2j}^*(P^*_j+P_j) }
= 0,\\
P_i(T,s) = Q_i(T,s),\\
\dot{\bar{P}}_i + Q_i+\bar{Q}_i + P (B_1 + \bar{B}_1)+(B_1 + \bar{B}_1)^*P_i
  \\
+ \sum_{s'\neq s}({\bar{P}}_i(t,s')-\bar{P}_i(t,s))\tilde{q}_{ss'} \\
- \frac{1}{4} (\bar{P}^*_i+\bar{P}_i)(B_{2i}+\bar{B}_{2i})(R_i+\bar{R}_i)^{-\frac{1}{2}*}(R_i+\bar{R}_i)^{-\frac{1}{2}}\\
(B_{2i}+\bar{B}_{2i})^* (\bar{P}^*_i+\bar{P}_i) \\
 -\frac{1}{4} \sum_{j \ne i} (\bar{P}^*_j+\bar{P}_j)^* (B_{2j}+\bar{B}_{2j}) (R_j+\bar{R}_j)^{-1*}\\
 (B_{2j}+\bar{B}_{2j})^* (\bar{P}^*_i+\bar{P}_i) \\
-\frac{1}{4} \sum_{j \ne i}  (\bar{P}^*_i+\bar{P}_i)(B_{2j}+\bar{B}_{2j}) (R_j+\bar{R}_j)^{-1*}\\
 (B_{2j}+\bar{B}_{2j})^*(\bar{P}^*_j+\bar{P}_j) 
= 0,\\
\bar{P}_i(T,s) = Q_i(T,s)+\bar{Q}_i(T,s),\\
\dot{\delta}_i + \frac{1}{2} \langle (P^*_i+P_i)S_0,S_0\rangle + \frac{1}{2} \int_{\Theta} \langle (P^*_i+P_i)M_0,M_0 \nu(d\theta) \rangle \\
+ \sum_{s'\neq s}({\delta}_i(t,s')-\delta_i(t,s))\tilde{q}_{ss'} = 0,\\
\delta_i(T,s) = 0,\\
\forall s \in \mathcal{S},
\end{array}
\right.
\end{equation}
whenever these differential equations have a unique  solution that does not blow up within $[0,T]$. \eod
\end{theorem}
Under the symmetric matrix assumption above, it is easy to check that if $P$ is a solution then $P^*$ is also a solution. Therefore $P^*_i(t,s)=P_i(t,s),\  (t,s)\in [0,T]\times \mathcal{S}.$

From the state system (\ref{equation:state}), the conditional expected matrix $\bar{X}(t):=\mathbb{E}[{X}(t)| \mathcal{F}^s_t],$ where  $\mathcal{F}^s$ is the natural filtration of the regime switching process $s$ up to $t,$ solves the following system:

\begin{align*}
\begin{array}{l}
d\bar{X} = [ (B_1 + \bar{B}_1)\bar{X} + \sum_{j=1}^{n} (B_{2j}+\bar{B}_{2j}) \bar{U}_j ]dt,\\
\bar{X}(0)=\mathbb{E}[X_0].
\end{array}
\end{align*}
which means that
\begin{align*}
\begin{array}{l}
d\bar{X} = [ (B_1 + \bar{B}_1) -  \frac{1}{2}\sum_{j=1}^{n} (B_{2j}+\bar{B}_{2j}) (R_j+\bar{R}_j)^{-1} \\
(B_{2j}+\bar{B}_{2j})^* (\bar{P}^*_j+\bar{P}_j) ] \bar{X}dt,\\
\bar{X}(0)=\mathbb{E}[X_0],
\end{array}
\end{align*} 
which will be used for feedback in the optimal strategy.  
Next, we provide a semi-explicit solution to the  full-cooperation case. 
\begin{coro} \label{cor:1}
 Assume that  $Q_0:=\sum_{j=1}^n Q_j,  Q_0+\bar{Q}_0:=\sum_{j=1}^n[Q_j+\bar{Q}_j], R_i,R_i+\bar{R}_i$ are symmetric positive definite. 
 The fully cooperative solution of the problem $\inf_{(U_1, \ldots, U_n)}\mathbb{E}[\sum_{j=1}^n L_j]$ is given by 
 
 \begin{align*}
\left\{
\begin{array}{l}
U_i^\mathrm{rn,g} - \bar{U}_i^\mathrm{rn,g} =- R_i^{-1}B_{2i}^* P_0(X - \bar{X}),\\
\bar{U}_i^\mathrm{rn,g} =  -(R_i+\bar{R}_i)^{-1} (B_{2i}+\bar{B}_{2i})^* \bar{P}_0\bar{X},\\
L_0^\mathrm{rn,g}(X,s,U^\mathrm{rn,g}) = \mathbb{E}\langle P_0(0,s(0)) (X_0 - \bar{X}_0), X_0 - \bar{X}_0  \rangle \\
+ \mathbb{E}\langle \bar{P}_0(0,s(0))  \bar{X}_0, \bar{X}_0  \rangle  + \mathbb{E}[\delta_0(0,s(0))],\\
i \in \{1,\dots, n\},
\end{array}
\right.
\end{align*}
where $P_0, \bar{P}_0$, and $\delta_0$ solve the following differential equations:
\begin{equation} \label{tt0}
\left\{
\begin{array}{l}
\dot{P}_0 + Q_0 +  P_0 B_1+B_1^*P_0 + \sum_{s'\neq s}({P}_0(t,s')-{P}_0(t,s))\tilde{q}_{ss'} \\
- P_0 \Big[ \sum_{\textcolor{black}{i=1}}^n B_{2i} R_i^{-1} B_{2i}^* \Big]  P_0= 0,\\
P_0(T,s) = Q_0(T,s),\\
\dot{\bar{P}}_0 + Q_0+\bar{Q}_0 + \bar{P}_0(B_1 + \bar{B}_1) + (B_1 + \bar{B}_1)^*\bar{P}_0 \\
+ \sum_{s'\neq s}({\bar{P}}_0(t,s')-\bar{P}_0(t,s))\tilde{q}_{ss'} \\
-\bar{P}_0\Big[\sum_{\textcolor{black}{i=1}}^n (B_{2i}+\bar{B}_{2i})(R_i+\bar{R}_i)^{-1}(B_{2i}+\bar{B}_{2i})^*\Big] \bar{P}_0  = 0,\\
\bar{P}_0(T,s) = Q_0(T,s)+\bar{Q}_0(T,s),\\
\dot{\delta}_0 + \langle P_0 S_0,S_0\rangle + \int_{\Theta} \langle P_0 M_0,M_0 \nu(d\theta) \rangle \\
+ \sum_{s'\neq s}({\delta}_0(t,s')-\delta_0(t,s))\tilde{q}_{ss'} = 0,\\
\delta_0(T,s) = 0,\\
\forall s \in \mathcal{S}.
\end{array}
\right.
\end{equation}
Notice that these Riccati equations have positive solution $P_0, \bar{P}_0$, and $\delta_0$ and there is no blow up in $[0,T]$. \eod
\end{coro}
The proof of Corollary \ref{cor:1} is immediate from Theorem \ref{thm:1} by  one single team and with a choice vector of matrices $U=(U_i)_{i\in \mathcal{I}}.$

\begin{coro} \label{cor:11}
 Assume that  $Q,  Q+\bar{Q}, R_i,R_i+\bar{R}_i$ are symmetric positive definite for $i\in\mathcal{I}_{+}$ and $ -R_j,-(R_j+\bar{R}_j)$ are symmetric positive definite for $j\in\mathcal{I}_{-}.$ We assume that 
 $\mathcal{I}_{+}\cup \mathcal{I}_{-}=\mathcal{I}.$
 The adversarial game problem of the team attackers $\mathcal{I}_{-}$ and the team of defenders  $\mathcal{I}_{+}$ has a saddle and it  is given by 
 
 \begin{align*}
\left\{
\begin{array}{l}
U_i^{\mathrm{ad}} - \bar{U}_i^{\mathrm{ad}} =- R_i^{-1}B_{2i}^* P^\mathrm{ad}(X - \bar{X}),\\ 
\bar{U}_i^\mathrm{ad} =  -(R_i+\bar{R}_i)^{-1} (B_{2i}+\bar{B}_{2i})^* \bar{P}^\mathrm{ad}\bar{X},\ i \in \mathcal{I}_{+}, \\
V_j^{\mathrm{ad}} - \bar{V}_j^{\mathrm{ad}} =- R_j^{-1}B_{2j}^* P^\mathrm{ad}(X - \bar{X}),\ j \in \mathcal{I}_{-},\\
 \bar{V}_j^\mathrm{ad} =  -(R_j+\bar{R}_j)^{-1} (B_{2j}+\bar{B}_{2j})^* \bar{P}^\mathrm{ad}\bar{X},\ j \in  \mathcal{I}_{-},\\
L^\mathrm{ad}(X,s,U^\mathrm{ad}) = \mathbb{E}\langle P^\mathrm{ad}(0,s(0)) (X_0 - \bar{X}_0), X_0 - \bar{X}_0  \rangle \\
+ \mathbb{E}\langle \bar{P}^\mathrm{ad}(0,s(0))  \bar{X}_0, \bar{X}_0  \rangle  + \mathbb{E}[\delta^\mathrm{ad}(0,s(0))],
\end{array}
\right.
\end{align*}
where $P^\mathrm{ad}, \bar{P}^\mathrm{ad}$, and $\delta^\mathrm{ad}$ solve the following differential equations:
\begin{equation} \label{tt0t}
\left\{
\begin{array}{l}
\dot{P} + Q+  P B_1+B_1^*P + \sum_{s'\neq s}({P}(t,s')-{P}(t,s))\tilde{q}_{ss'} \\
- P {\color{blue}\Big[ \sum_{i\in  \mathcal{I}_{+}} B_{2i} R_i^{-1} B_{2i}^*+ \sum_{j\in \mathcal{I}_{-}}B_{2j} R_j^{-1} B_{2j}^* \Big] } P= 0,\\
P(T,s) = Q(T,s),\\
\dot{\bar{P}} + Q+\bar{Q} + \bar{P}(B_1 + \bar{B}_1) + (B_1 + \bar{B}_1)^*\bar{P}\\
+ \sum_{s'\neq s}({\bar{P}}(t,s')-\bar{P}(t,s))\tilde{q}_{ss'} \\
-\bar{P}{\color{blue} \Big[\sum_{i\in \mathcal{I}_{+}} (B_{2i}+\bar{B}_{2i})(R_i+\bar{R}_i)^{-1}(B_{2i}+\bar{B}_{2i})^* }\\
{\color{blue}+ \sum_{j\in \mathcal{I}_{-}} (B_{2j}+\bar{B}_{2j})(R_j+\bar{R}_j)^{-1}(B_{2j}+\bar{B}_{2j})^*\Big]} \bar{P}  = 0,\\
\bar{P}(T,s) = Q(T,s)+\bar{Q}(T,s),\\
\dot{\delta} + \langle P S_0,S_0\rangle + \int_{\Theta} \langle P M_0,M_0 \nu(d\theta) \rangle \\
+ \sum_{s'\neq s}({\delta}(t,s')-\delta(t,s))\tilde{q}_{ss'} = 0,\\
\delta(T,s) = 0,\\
\forall s \in \mathcal{S}.
\end{array}
\right.
\end{equation}\eod
\end{coro}
The proof of Corollary \ref{cor:11} is immediate from Theorem \ref{thm:1} by  considering two adversarial teams and with  choice vector of matrices $U_+=(U_i)_{i\in \mathcal{I}_+}$ and $U_-=(U_i)_{i\in \mathcal{I}_-}$ respectively.

Notice that the Riccati equations in (\ref{tt0t}) have positive definite solution $P^\mathrm{ad}$ if in addition 
\begin{align*}
\begin{array}{l}
{\color{blue}\Big[ \sum_{i\in  \mathcal{I}_{+}} B_{2i} R_i^{-1} B_{2i}^*+ \sum_{j\in \mathcal{I}_{-}}B_{2j} R_j^{-1} B_{2j}^* \Big] \succ 0},
\end{array}
\end{align*}
which does not  blow up within $[0,T]$, and positive solution $\bar{P}^\mathrm{ad}$ if in addition

\begin{align*}
\begin{array}{l}
\Big[\sum_{i\in \mathcal{I}_{+}} (B_{2i}+\bar{B}_{2i})(R_i+\bar{R}_i)^{-1}(B_{2i}+\bar{B}_{2i})^* \\
+ \sum_{j\in \mathcal{I}_{-}} (B_{2j}+\bar{B}_{2j})(R_j+\bar{R}_j)^{-1}(B_{2j}+\bar{B}_{2j})^*\Big]\succ 0,
\end{array}
\end{align*}
within $[0,T].$
Next, we study the risk-sensitive case and point out some facts regarding the comparison of its solution with respect to the risk-neutral case as the risk-sensitivity index vanishes.

\subsection{Risk-Sensitive Case}
A risk averse decision-maker (with cost functional) is a decision-maker who prefers higher cost with known risks rather than lower cost with unknown risks. 
In other words, among various control strategies giving the same cost  with different level of risks, this decision-maker always prefers the alternative with the lowest risk. 

When $M_0\neq 0$ the exponential martingale of compensated Poisson random process times a linear process yields to an exponential non-quadratic terms
$\mathbb{E}\mathrm{exp}[{\int_0^T \int_{\Theta}\langle (P^*_i +P_i) (X-\bar{X}), M_0\tilde{N}(dt,d\theta)\rangle}]$  has an  exponential non-quadratic term. 
Therefore we consider the risk-sensitive case when $M_0$ vanishes (no Poisson jump) and a single regime $\mathcal{S}=\{s_0\}.$  When $\lambda_i>0$ the decision-maker $i$ is risk-averse and when $\lambda_i<0$ is risk-seeking. As $\lambda_i$ goes to zero, decision-maker $i$ becomes a risk-neutral decision-maker.  The best response problem of decision-maker $i$  is well-posed only for $\lambda_i \leq \bar{\lambda}_i$ where $\bar{\lambda}_i$ will be determined from the solution region of the differential system derived below.

\begin{theorem} \label{thm:2}
Assume that  $Q_i, R_i, Q_i+\bar{Q}_i, R_i+\bar{R}_i$ are symmetric positive definite. There is no jump $M_0=0$  and a single regime $\mathcal{S}=\{s_0\}.$ 
Then the matrix-valued mean-field-type (risk-sensitive) Nash equilibrium strategy and the (risk-sensitive) equilibrium  cost are given by:
\begin{align*}
\left\{
\begin{array}{l}
U^\mathrm{rs}_i - \bar{U}^\mathrm{rs}_i =-  R_i^{-1}B_{2i}^*P^\mathrm{rs}_i(X - \bar{X}),\\
\bar{U}^\mathrm{rs}_i=  - (R_i+\bar{R}_i)^{-1} (B_{2i}+\bar{B}_{2i})^* \bar{P}_i\bar{X},\\
L_i^\mathrm{rs}(X,s,U^\mathrm{rs}) \\
= \frac{1}{\lambda_i}\log \mathbb{E} \exp\{\lambda_i [\langle P^\mathrm{rs}_i(0,s(0)) (X_0 - \bar{X}_0), X_0 - \bar{X}_0  \rangle \\
+\langle \bar{P}_i(0,s(0))  \bar{X}_0, \bar{X}_0  \rangle  + \delta^\mathrm{rs}_i(0,s(0))]\},\\
 i \in \{1,\dots, n\},
\end{array}
\right.
\end{align*}
where $P^\mathrm{rs}_i, \bar{P}^\mathrm{rs}_i= \bar{P}_i$, and $\delta^\mathrm{rs}_i$ solve the following differential equations:
\begin{equation} \label{tt}
\left\{
\begin{array}{l}
\dot{P}_i + Q_i +P_i  B_1^*+ B_1P_i \\ 
- P_i {\color{blue} [ B_{2i}  R_i^{-1} B_{2i}^* -2\lambda_i S_0S_0^*  ] }P_i \\
-  \sum_{j \ne i} P_j B_{2j} R_j^{-1*}B_{2j}^*P_i 
-  \sum_{j \ne i}  P_iB_{2j} R_j^{-1*}B_{2j}^* P_j= 0,\\
P_i(T,s) = Q_i(T,s),\\
\dot{\bar{P}}_i + Q_i+\bar{Q}_i + \bar{P}_i(B_1 + \bar{B}_1)^*+ (B_1 + \bar{B}_1)\bar{P}_i \\
- \bar{P}_i{\color{black}[ (B_{2i}+\bar{B}_{2i})(R_i+\bar{R}_i)^{-1}(B_{2i}+\bar{B}_{2i})^* ] }\bar{P}_i \\
-  \sum_{j \ne i} \bar{P}_j (B_{2j}+\bar{B}_{2j}) (R_j+\bar{R}_j)^{-1*} (B_{2j}+\bar{B}_{2j})^* \bar{P}_i\\
-  \sum_{j \ne i}  \bar{P}_i(B_{2j}+\bar{B}_{2j}) (R_j+\bar{R}_j)^{-1*} (B_{2j}+\bar{B}_{2j})^*\bar{P}_j 
= 0,\\
\bar{P}_i(T,s) = Q_i(T,s)+\bar{Q}_i(T,s),\\
\dot{\delta}_i + \langle P_i S_0,S_0\rangle =0, 
\delta_i(T,s) = 0,\\
\forall s \in \mathcal{S},
\end{array}
\right.
\end{equation}
whenever these differential system of equations have a unique  solution that does not blow up in $[0,T]$. \eod
\end{theorem}

As all  $\lambda_i$ vanish, the matrix-valued differential system  (\ref{tt}) becomes  the risk-neutral system (\ref{tt0}) above, and the risk-sensitive optimal strategy coincides with the risk-neutral one. A bound for $\bar{\lambda}_i$ can be obtained from the positivity condition of the matrices   
\begin{align}
\label{condition}
[ B_{2i}  R_i^{-1} B_{2i}^* -2\lambda_i S_0S_0^*  ] \succ 0.
\end{align}
Note that $\bar{P}^\mathrm{rs}_i= \bar{P}_i$ because this coefficient is associated with the  term $\langle \bar{X},  \bar{X}\rangle$ which is independent of the Brownian motion.

Consider a  finite population of decision-makers $\mathcal{I}:= \{1,\ldots n\} $ is classified as follows:
\begin{itemize}
\item  risk-neutral decision-makers: $\mathcal{I}_{0}=\{ i\in \{1,\ldots n\}\  | \ \lambda_i\rightarrow 0\}$
 \item risk-averse  decision-makers: $\mathcal{I}_{+}=\{ i\ | \ \lambda_i>0\}$
 \item risk-seeking decision-makers: $\mathcal{I}_{-}=\{ i\ | \  \lambda_i<0\}$
\end{itemize}
\begin{coro} \label{cor:22} A mixture of risk-neutral, risk-seeking and risk-averse are obtained solving the following system:

\begin{align*}
&\left\{\begin{array}{l}
i\in \mathcal{I}_0:\\
\dot{P}_i + Q_i +P_i  B_1^*+ B_1P_i\\  
- P_i   B_{2i}  R_i^{-1} B_{2i}^*  P_i -  \sum_{j \ne i} P_j B_{2j} R_j^{-1*}B_{2j}^*P_i \\
- \sum_{j \ne i} P_i B_{2j} R_j^{-1*}B_{2j}^*P_j = 0,\\ 
\end{array}\right.\\
&\left\{\begin{array}{l}
i\in \mathcal{I}_+:\\
\dot{P}_i + Q_i +P_i  B_1^*+ B_1P_i\\  
- P_i {\color{blue} [ B_{2i}  R_i^{-1} B_{2i}^* -2\lambda_i S_0S_0^*  ] }P_i \\
-  \sum_{j \ne i} P_j B_{2j} R_j^{-1*}B_{2j}^*P_i\\
-  \sum_{j \ne i} P_i B_{2j} R_j^{-1*}B_{2j}^*P_j = 0,\\ 
\end{array}\right.\\
&\left\{\begin{array}{l}
i\in \mathcal{I}_{-}:\\
\dot{P}_i + Q_i +P_i  B_1^*+ B_1P_i\\ 
- P_i {\color{blue} [ B_{2i}  R_i^{-1} B_{2i}^* +2(-\lambda_i) S_0S_0^*  ] }P_i \\
-  \sum_{j \ne i} P_j B_{2j} R_j^{-1*}B_{2j}^*P_i\\
-  \sum_{j \ne i} P_i B_{2j} R_j^{-1*}B_{2j}^*P_j = 0,\\
\end{array}\right.\\
&P_i(T,s) = Q_i(T,s), i\in \mathcal{I},\\
&\forall s \in \mathcal{S}.
\end{align*}
\end{coro}

\begin{coro} \label{cor:2}
Assume that  $Q_0:=\sum_{j=1}^n Q_i,  Q_0+\bar{Q}_0:=\sum_{j=1}^n[Q_i,+\bar{Q}_i], R_i,R_i+\bar{R}_i$ are symmetric positive definite.  The risk-sensitive fully cooperative solution of the  problem $\inf_{(U_1, \ldots, U_n)}\frac{1}{\lambda}\log \mathbb{E}[e^{\lambda\sum_{j=1}^n L_j}]$ is given by 

 \begin{align*}
\left\{
\begin{array}{l}
U_i^\mathrm{rs,g} - \bar{U}_i^\mathrm{rs,g} =- R_i^{-1}B_{2i}^* P^\mathrm{rs,g}_0(X - \bar{X}),\\
\bar{U}_i^\mathrm{rs,g} =  -(R_i+\bar{R}_i)^{-1} (B_{2i}+\bar{B}_{2i})^* \bar{P}_0\bar{X},\\
L_0^\mathrm{rs,g}(X,s,U^\mathrm{rs,g}) \\
= \frac{1}{\lambda}\log \mathbb{E} \exp\{\lambda [\langle P^\mathrm{rs,g}_0(0,s(0)) (X_0 - \bar{X}_0), X_0 - \bar{X}_0  \rangle \\
+\langle \bar{P}_0(0,s(0))  \bar{X}_0, \bar{X}_0  \rangle  + \delta^\mathrm{rs,g}_0(0,s(0))]\},\\
 i\in \mathcal{I}, 
\end{array}
\right.
\end{align*}
where $P^\mathrm{rs,g}_0, \bar{P}^\mathrm{rs,g}_0=\bar{P}_0$, and $\delta^\mathrm{rs,g}_0$ solve the following differential equations:
\begin{equation} \label{tt0}
\left\{
\begin{array}{l}
\dot{P}_0 + Q_0 +  P_0 B_1+B_1^*P_0 \\ 
- P_0 \textcolor{blue}{ \Big[ \sum_{i=1}^n B_{2i} R_i^{-1} B_{2i}^* - 2 \lambda S_0S_0^*  \Big]}  P_0= 0,\\
P_0(T,s) = Q_0(T,s),\\
\dot{\bar{P}}_0 + Q_0+\bar{Q}_0 + \bar{P}_0(B_1 + \bar{B}_1) + (B_1 + \bar{B}_1)^*\bar{P}_0 \\
-\bar{P}_0\Big[\sum_{i=1}^n (B_{2i}+\bar{B}_{2i})(R_i+\bar{R}_i)^{-1}(B_{2i}+\bar{B}_{2i})^*\Big] \bar{P}_0  = 0,\\
\bar{P}_0(T,s) = Q_0(T,s)+\bar{Q}_0(T,s),\\
\dot{\delta}_0 + \langle P_0 S_0,S_0\rangle =0, 
\delta_0(T,s) = 0,\\
\forall s \in \mathcal{S}.
\end{array}
\right.
\end{equation}
 \eod

\end{coro}

\begin{remark} (Shared-Risk Situation)  
Notice that the bound $\bar{\lambda}$ is obtained from the positivity condition $\sum_{i=1}^n B_{2i} R_i^{-1} B_{2i}^* - 2 \lambda S_0S_0^*  \succ 0$. Hence, the risk condition is relaxed thanks to the full cooperation in comparison with the non-cooperative risk consideration, i.e.,
\begin{align*}
\bar{\lambda} &= \mathrm{sup} \left\{ \lambda | \begin{array}{c}\text{full cooperation risk-sensitive}\\ \text{problem is well-posed} \end{array} \right\},\\
\bar{\lambda}_i &= \mathrm{sup} \left\{ \lambda_i | \begin{array}{c}\text{non-cooperative risk-sensitive}\\ \text{problem is well-posed} \end{array} \right\},
\end{align*}
where 
\begin{align*}
\begin{array}{l}
B_{2i} R_i^{-1} B_{2i}^* \succ 2 \lambda_i S_0S_0,\;\; (\text{from \eqref{condition}}),\\
\sum_{i=1}^n B_{2i} R_i^{-1} B_{2i}^* \succ 2 \sum_{i=1}^n\lambda_i S_0S_0,
\end{array}
\end{align*}
therefore, it is concluded that $\sum_{i=1}^n\bar{\lambda}_i \leq \bar{\lambda}.$ Full cooperation increases the well-posedness domain by means shared risk. \eod
\end{remark}
  
\begin{coro} \label{cor:11rsad}
 Assume that  $Q,  Q+\bar{Q}, R_i,R_i+\bar{R}_i$ are symmetric positive definite for $i\in\mathcal{I}_{+}$ and $ -R_j,-(R_j+\bar{R}_j)$ are symmetric positive definite for $j\in\mathcal{I}_{-}.$ We assume that 
 $\mathcal{I}_{+}\cup \mathcal{I}_{-}=\mathcal{I}.$
 The adversarial risk-sensitive game problem of the team attackers $\mathcal{I}_{-}$ and the team of defenders  $\mathcal{I}_{+}$ has a saddle and it  is given by 
 
 \begin{align*}
\left\{
\begin{array}{l}
U_i^{\mathrm{rs,ad}} - \bar{U}_i^{\mathrm{rs,ad}} =- R_i^{-1}B_{2i}^* P^{\mathrm{rs,ad}}(X - \bar{X}),\\ 
\bar{U}_i^\mathrm{rs,ad} =  -(R_i+\bar{R}_i)^{-1} (B_{2i}+\bar{B}_{2i})^* \bar{P}^{\mathrm{rs,ad}}\bar{X},\ i \in \mathcal{I}_{+}, \\
V_j^{\mathrm{rs,ad}} - \bar{V}_j^{\mathrm{rs,ad}} =- R_j^{-1}B_{2j}^* P^{\mathrm{rs,ad}}(X - \bar{X}),\\ 
\bar{V}_j^\mathrm{rs,ad} =  -(R_j+\bar{R}_j)^{-1} (B_{2j}+\bar{B}_{2j})^* \bar{P}^{\mathrm{rs,ad}}\bar{X},\ j \in  \mathcal{I}_{-},\\
L^\mathrm{ad}(X,s,U^\mathrm{ad}) = \mathbb{E}\langle P^\mathrm{rs,ad}(0,s(0)) (X_0 - \bar{X}_0), X_0 - \bar{X}_0  \rangle \\
+ \mathbb{E}\langle \bar{P}^\mathrm{ad}(0,s(0))  \bar{X}_0, \bar{X}_0  \rangle  + \mathbb{E}[\delta^\mathrm{rs,ad}(0,s(0))],
\end{array}
\right.
\end{align*}
where $P^\mathrm{rs,ad}, \bar{P}^\mathrm{rs,ad}=\bar{P}^\mathrm{ad}$, and $\delta^{\mathrm{rs,ad}}$ solve the following differential equations:
\begin{equation} \label{tt0t}
\left\{
\begin{array}{l}
\dot{P} + Q+  P B_1+B_1^*P + \sum_{s'\neq s}({P}(t,s')-{P}(t,s))\tilde{q}_{ss'} \\
- P {\color{blue}\Big[ \sum_{i\in  \mathcal{I}_{+}} B_{2i} R_i^{-1} B_{2i}^*} \\
{\color{blue}+ \sum_{j\in \mathcal{I}_{-}}B_{2j} R_j^{-1} B_{2j}^* -2 \lambda S_0S_0^* \Big] } P= 0,\\
P(T,s) = Q(T,s),\\
\dot{\bar{P}} + Q+\bar{Q} + \bar{P}(B_1 + \bar{B}_1) + (B_1 + \bar{B}_1)^*\bar{P} \\
+ \sum_{s'\neq s}({\bar{P}}(t,s')-\bar{P}(t,s))\tilde{q}_{ss'} \\
-\bar{P}{\color{blue} \Big[\sum_{i\in \mathcal{I}_{+}} (B_{2i}+\bar{B}_{2i})(R_i+\bar{R}_i)^{-1}(B_{2i}+\bar{B}_{2i})^*}\\ 
{\color{blue}+ \sum_{j\in \mathcal{I}_{-}} (B_{2j}+\bar{B}_{2j})(R_j+\bar{R}_j)^{-1}(B_{2j}+\bar{B}_{2j})^*}\\
{\color{blue}-2 \lambda S_0S_0^*\Big]} \bar{P}  = 0,\\
\bar{P}(T,s) = Q(T,s)+\bar{Q}(T,s),\\
\dot{\delta} + \langle P S_0,S_0\rangle + \sum_{s'\neq s}({\delta}(t,s')-\delta(t,s))\tilde{q}_{ss'} = 0,\\
\delta(T,s) = 0,\\
\forall s \in \mathcal{S}.
\end{array}
\right.
\end{equation}\eod
\end{coro}
Notice that these risk-sensitive adversarial Riccati equations have positive definite solution $P^{\mathrm{rs,ad}}$ if in addition 
\begin{align*}
\begin{array}{l}
{\color{blue}\Big[ \sum_{i\in  \mathcal{I}_{+}} B_{2i} R_i^{-1} B_{2i}^*+ \sum_{j\in \mathcal{I}_{-}}B_{2j} R_j^{-1} B_{2j}^* -2 \lambda S_0S_0^*\Big] \succ 0},
\end{array}
\end{align*}
which does not  blow up within $[0,T].$
 
\section{Relationship Between Adversarial and Risk-Sensitive Game Problem}
We discuss the relationship between the solutions obtained from an adversarial game problem and the risk-sensitive game problem when there is no jump and a single regime.

\subsection{Mean-Field-Free Control}

Consider the following risk-sensitive control problem with $\lambda>0$:
\begin{align*}
\begin{array}{l}
L^{\mathrm{rs}}(X^\mathrm{rs},s,U^\mathrm{rs})=\frac{1}{\lambda}\log \left[ \mathbb{E}e^{\lambda L(X^\mathrm{rs},s,U^\mathrm{rs})} \right],\\
\inf_{U^\mathrm{rs}}  L^{\mathrm{rs}}(X^\mathrm{rs},s,U^\mathrm{rs})\\
\text{s.t.}~~dX^\mathrm{rs} = [B_1X^\mathrm{rs}+B_2U^\mathrm{rs}]dt + S_0dB,\\ X(0)=X_0,\ s(0)=s_0
\end{array}
\end{align*}
where
\begin{align*}
\begin{array}{l}
L(X^\mathrm{rs},s,U^\mathrm{rs}) = \langle Q(T,s(T)) {X}^\mathrm{rs}(T), {X}^\mathrm{rs}(T)  \rangle \\
+ \int_{0}^{T} \langle Q X^\mathrm{rs}, X^\mathrm{rs} \rangle  dt + \int_{0}^{T} \langle R U^\mathrm{rs}, U^\mathrm{rs}  \rangle  dt,
\end{array}
\end{align*}
Its solution is given by:
\begin{align*}
\begin{array}{l}
U^\mathrm{rs} = -R^{-1}B_2^*P^\mathrm{rs} X^\mathrm{rs},\\ 
L^\mathrm{rs}(X^\mathrm{rs},s,U^\mathrm{rs}) \\
= \frac{1}{\lambda}\log \mathbb{E} \exp\{\lambda [\langle P^\mathrm{rs}_0(0,s(0)) X^\mathrm{rs}_0, X^\mathrm{rs}_0  \rangle  + \delta^\mathrm{rs}_0(0,s(0))]\},\\
\end{array}
\end{align*}
Now  consider the following adversarial risk-neutral control problem:
\begin{align*}
\begin{array}{l}
\inf_{U^\mathrm{ad}} \sup_{V^\mathrm{ad}}\mathbb{E}[L^{\mathrm{ad}}(U,V)],\\
\text{s. t.}~~ dX^\mathrm{ad} = \Big[B_1X^\mathrm{ad}+B_2U^\mathrm{ad}+\Big[\sqrt{2\lambda}S_0[-\tilde{R}]^{\frac{1}{2}}\Big]V^\mathrm{ad}\Big]dt \\
+ S_0dB,
\end{array}
\end{align*}
with
\begin{align*}
\begin{array}{l}
L^\mathrm{ad}(X^\mathrm{ad},s,U^\mathrm{ad}) =  L(X^\mathrm{ad},s,U^\mathrm{ad}) 
+ \int_{0}^{T} \langle \tilde{R} V^\mathrm{ad}, V^\mathrm{ad}  \rangle  dt,
\end{array}
\end{align*}
where  $\tilde{R}\prec 0$. Its explicit solution is given by
\begin{align*}
\begin{array}{l}
U^\mathrm{ad} = -R^{-1}B_2^*P^\mathrm{ad}X^\mathrm{ad},\\
V^\mathrm{ad} = \sqrt{2\lambda}[-\tilde{R}]^{-\frac{1}{2}}S_0^*P^\mathrm{ad}X^\mathrm{ad},\\
L^\mathrm{ad}(X^\mathrm{ad},s,U^\mathrm{ad}) = \mathbb{E}\langle P^\mathrm{ad}(0,s(0))X^\mathrm{ad}_0,X^\mathrm{ad}_0 \rangle \\
+ \mathbb{E}[\delta^\mathrm{ad}(0,s(0))] \\
\end{array}
\end{align*}
the coefficients $P^\mathrm{ad}$ and $ P^\mathrm{rs}$ solve the same differential system which is given by
\begin{align*}
\begin{array}{l}
\dot{P} + Q +  P B_1+B_1^*P - P \Big[ B_{2} R^{-1} B_{2}^* - 2 \lambda S_0S_0^*  \Big]  P= 0,\\ P(T, s)=Q(T, s),
\end{array}
\end{align*}

{\color{blue} However, the evolution of the system states $X^\mathrm{rs}$ and $X^\mathrm{ad}$ may not have the same probability law even if they start at the same initial random variable $X_0.$
When $X_0$ is a deterministic matrix and $\mathcal{S}$ is a singleton, the two optimal costs  $L^\mathrm{ad}$ and $L^\mathrm{rs}$ coincide. In general, these costs are different.

By taking the expectation, it follows that
   
   $d\bar{X}^\mathrm{rs} = [B_1-B_{2}R^{-1}B_{2}^*P] \bar{X}^\mathrm{rs}dt,$ and

   $$d\bar{X}^\mathrm{ad} = [B_1-(B_{2}R^{-1}B_{2}^*-2\lambda S_0S_0^*)P] \bar{X}^\mathrm{ad}dt$$ 
   
Hence $$\bar{X}^\mathrm{ad}(t)=\bar{X}^\mathrm{rs}(t) e^{2\lambda \int_0^t S_0S_0^*P(t') dt'},$$
 We conclude that the first moments of the two problems are different for $\lambda S_0S_0^*\neq 0 $ and $\mathbb{E}X_0\neq 0.$
 This means that, in a problem the first moment is involved in the cost, the mean of the mean-field trajectory of risk-sensitive and robust design are not necessary equal. Below we show 
 that one can modify the robust control law such that the first moment of these designs are equal. It will be achieved by means   of mean-field-type robust control.  }

\subsection{Mean-Field-Free Team}

Consider the following risk-sensitive team problem with team members in $\mathcal{I}_+$ :
\begin{align*}
\begin{array}{l}
\inf_{U^\mathrm{rs}} \frac{1}{\lambda}\log \left[ \mathbb{E}e^{\lambda L_{\mathrm{team}}(X^\mathrm{rs},s,U^\mathrm{rs})} \right],\\
\text{s. t.}~~
dX^\mathrm{rs} = [B_1X^\mathrm{rs}+\sum_{i \in \mathcal{I}_+}B_{2i}U_i^\mathrm{rs}]dt + S_0dB,
\end{array}
\end{align*}
where
\begin{align*}
\begin{array}{l}
L_{\mathrm{team}}(X^\mathrm{rs},s,U^\mathrm{rs}) = \langle Q(T,s(T)) {X}^\mathrm{rs}(T), {X}^\mathrm{rs}(T)  \rangle \\
+ \int_{0}^{T} \langle Q X^\mathrm{rs}, X^\mathrm{rs} \rangle  dt + \int_{0}^{T} \sum_{i \in \mathcal{I}_+} \langle R_i U_i^\mathrm{rs}, U_i^\mathrm{rs}  \rangle  dt,
\end{array}
\end{align*}
Its solution is given by:
\begin{align*}
\begin{array}{l}
U_i^\mathrm{rs} = -R^{-1}_iB_{2i}^*P^\mathrm{rs} X^\mathrm{rs},\\ 
L^\mathrm{rs}_{\mathrm{team}}(X^\mathrm{rs},s,U^\mathrm{rs}) \\
= \frac{1}{\lambda}\log \mathbb{E} \exp\{\lambda [\langle P^\mathrm{rs}_0(0,s(0)) X^\mathrm{rs}_0, X^\mathrm{rs}_0  \rangle  + \delta^\mathrm{rs}_0(0,s(0))]\}.\\
\end{array}
\end{align*}

We now design  the following adversarial risk-neutral game problem with two teams in $I_{+}$ and $I_{-}$  respectively 
\begin{align*}
\begin{array}{l}
\inf_{U^\mathrm{ad}} \sup_{V^\mathrm{ad}}\mathbb{E}[L^{\mathrm{ad}}_{\mathrm{team}}(U,V)],\\
\text{s. t.}~~dX^\mathrm{ad} = \Big[B_1X^\mathrm{ad}+\sum_{i \in \mathcal{I}_+}B_{2i}U^\mathrm{ad} \\
+ |\mathcal{I}_-|^{-1}\sum_{j \in \mathcal{I}_-}\Big[\sqrt{2\lambda}S_0[-\tilde{R}_j]^{\frac{1}{2}}\Big]V_j^\mathrm{ad}\Big]dt + S_0dB, 
\end{array}
\end{align*}
with
\begin{align*}
\begin{array}{l}
L^\mathrm{ad}_{\mathrm{team}}(X^\mathrm{ad},s,U^\mathrm{ad}) = \langle Q(T,s(T)) {X}^\mathrm{ad}(T), {X}^\mathrm{ad}(T)  \rangle \\
+ \int_{0}^{T} \langle Q X^\mathrm{ad}, X^\mathrm{ad} \rangle  dt 
+ \int_{0}^{T} \sum_{i \in \mathcal{I}_+} \langle R_i U_i^\mathrm{ad}, U_i^\mathrm{ad}  \rangle  dt \\
+ \int_{0}^{T} \sum_{j \in \mathcal{I}_-} \langle \tilde{R}_j V_j^\mathrm{ad}, V_j^\mathrm{ad}  \rangle  dt,
\end{array}
\end{align*}
$\tilde{R}_j \prec 0, \ j\in \mathcal{I}_{-}$. Its explicit solution is given by
\begin{align*}
\begin{array}{l}
U_i^\mathrm{ad} = -R_i^{-1}B_{2i}^*P^\mathrm{ad}X^\mathrm{ad},\\
V_j^\mathrm{ad} = \sqrt{2\lambda}[-\tilde{R}]_j^{\frac{1}{2}}S_0^*P^\mathrm{ad}X^\mathrm{ad},\\
L^\mathrm{ad}_{\mathrm{team}}(X^\mathrm{ad},s,U^\mathrm{ad}) = \mathbb{E}\langle P^\mathrm{ad}(0,s(0))X^\mathrm{ad}_0,X^\mathrm{ad}_0 \rangle \\
+ \mathbb{E}[\delta^\mathrm{ad}(0,s(0))]. \\
\end{array}
\end{align*}
The coefficients $P^\mathrm{ad}$ and $ P^\mathrm{rs}$ solve the same differential equation which is given by
\begin{align*}
\begin{array}{l}
\dot{P} + Q +  P B_1+B_1^*P \\
- P  \Big[ \sum_{i \in \mathcal{I}_+} B_{2i} R_i^{-1} B_{2i}^* - 2 \lambda S_0S_0^*  \Big]  P= 0,\\ P(T,s)=Q(T,s),
\end{array}
\end{align*}
   The risk-sensitive optimal control  $U^\mathrm{rs}$ and the risk-neutral adversarial team optimal strategy $U^\mathrm{rs}$ are the similar structure. However the probability law of
   the system states $X^\mathrm{rs}$ and $X^\mathrm{ad}$ may be different even when they start from the same initial random variable $X_0.$ 
\subsection{Mean-Field-Type Team Problem}
We now examine possible  relationship between risk-sensitive team (full cooperation)  and robust team game in the  mean-field dependent case.  
Consider the following risk-sensitive mean-field-type team problem, involving the teams $\mathcal{I}_+$:

\begin{align*}
\begin{array}{l}
\inf_{U^\mathrm{rs}} \frac{1}{\lambda}\log \left[ \mathbb{E}e^{\lambda L_{\mathrm{team}}(X,s,U^\mathrm{rs})} \right],\\
\text{s. t.}~~dX = [ B_1(X - \bar{X}) + (B_1 + \bar{B}_1)\bar{X} \\
+ \sum_{i \in \mathcal{I}_+} B_{2i} (U_i - \bar{U}_i) + \sum_{i \in \mathcal{I}_+} (B_{2i}+\bar{B}_{2i}) \bar{U}_i ]dt\\
+ S_0 dB,
\end{array}
\end{align*}
where
\begin{align*}
\begin{array}{l}
L_{\mathrm{team}}(X,s,(U_i)_{i\in  \mathcal{I}_{+}} ) \\
= \langle Q(T,s(T)) ({X}(T)-\bar{X}(T)), {X}(T) -\bar{X}(T) \rangle \\
+ \langle (Q(T,s(T) + \bar{Q}(T,s(T)) \bar{X}(T), \bar{X}(T) \rangle\\
+ \int_{0}^{T} \langle Q (X-\bar{X}), X-\bar{X} \rangle  dt \\
+ \int_{0}^{T} \langle (Q + \bar{Q}) \bar{X}, \bar{X} \rangle  dt
+ \int_{0}^{T} \sum_{i\in \mathcal{I}_{+}}\langle R_i (U_i-\bar{U}_i), U_i-\bar{U}_i \rangle  dt \\
+ \int_{0}^{T} \sum_{i\in \mathcal{I}_{+}}\langle (R_i + \bar{R}_i) \bar{U}_i, \bar{U}_i \rangle  dt.
\end{array}
\end{align*}
Its solution is given by
 \begin{align*}
\left\{
\begin{array}{l}
U_i^\mathrm{rs} - \bar{U}_i^\mathrm{rs} =- R_i^{-1} B_{2i}^* P^\mathrm{rs}(X^\mathrm{rs} - \bar{X}^\mathrm{rs}),\\
\bar{U}_i^\mathrm{rs} =  -(R_i+\bar{R}_i)^{-1} (B_{2i}+\bar{B}_{2i})^* \bar{P}^\mathrm{rs}\bar{X}^\mathrm{rs},\\
L^\mathrm{rs}_{\mathrm{team}}(X^\mathrm{rs},s,U^\mathrm{rs}) \\
= \frac{1}{\lambda}\log \mathbb{E} \exp\{\lambda [\langle P^\mathrm{rs}(0,s(0)) (X^\mathrm{rs}_0 - \bar{X}^\mathrm{rs}_0), X^\mathrm{rs}_0 - \bar{X}^\mathrm{rs}_0  \rangle \\
+\langle \bar{P}^\mathrm{rs}(0,s(0))  \bar{X}^\mathrm{rs}_0, \bar{X}^\mathrm{rs}_0  \rangle  + \delta^\mathrm{rs}_0(0,s(0))]\},\\
i \in \mathcal{I}_+.
\end{array}
\right.
\end{align*}
We design the following adversarial risk-neutral mean-field-type game problem with two teams in $ \mathcal{I}_+$ and $ \mathcal{I}_{-}$:
\begin{align*}
\begin{array}{l}
\inf_{U^\mathrm{ad}} \sup_{V^\mathrm{ad}}\mathbb{E}[L^{\mathrm{ad}}_{\mathrm{team}}(U,V)],\\
\text{s. t.}~~dX^\mathrm{ad} = [ B_1(X^\mathrm{ad} - \bar{X}^\mathrm{ad}) + (B_1 + \bar{B}_1)\bar{X}^\mathrm{ad} \\
+ \sum_{i \in \mathcal{I}_+} B_{2i} (U_i^\mathrm{ad} - \bar{U}_i^\mathrm{ad}) + \sum_{i \in \mathcal{I}_+} (B_{2i}+\bar{B}_{2i}) \bar{U}_i^\mathrm{ad} \\
+ |\mathcal{I}_-|^{-1} \sum_{j \in \mathcal{I}_-}  \sqrt{2\lambda}S_0[-\tilde{R}_j]^{\frac{1}{2}} (V_j^\mathrm{ad} - \bar{V}_j^\mathrm{ad})  ]dt + S_0 dB ,
\end{array}
\end{align*}
with
\begin{align*}
\begin{array}{l}
L^\mathrm{ad}_{\mathrm{team}}(X^\mathrm{ad},s,U^\mathrm{ad}) = L_{\mathrm{team}}(X^\mathrm{ad},s,U^\mathrm{ad})\\
+ \int_{0}^{T} \sum_{j \in \mathcal{I}_-} \langle \tilde{R}_j (V_j^\mathrm{ad}-\bar{V}_j^\mathrm{ad}), V_j^\mathrm{ad}-\bar{V}_j^\mathrm{ad}  \rangle  dt,
\end{array}
\end{align*}
with $\tilde{R}_j \prec 0$. Its explicit solution is given by
 \begin{align*}
\left\{
\begin{array}{l}
U_i^{\mathrm{ad}} - \bar{U}_i^{\mathrm{ad}} =- R_i^{-1}B_{2i}^* P^\mathrm{ad}(X^\mathrm{ad} - \bar{X}^\mathrm{ad}),\\ 
\bar{U}_i^\mathrm{ad} =  -(R_i+\bar{R}_i)^{-1} (B_{2i}+\bar{B}_{2i})^* \bar{P}^\mathrm{ad}\bar{X}^\mathrm{ad},\ i \in \mathcal{I}_{+}, \\
V_j^{\mathrm{ad}} - \bar{V}_j^{\mathrm{ad}} = \sqrt{2\lambda}[-\tilde{R}_j]^{-\frac{1}{2}}S_0^*P^\mathrm{ad}(X^\mathrm{ad}-\bar{X}^\mathrm{ad}), \ j \in \mathcal{I}_{-},\\
\bar{V}_j^{\mathrm{ad}}=0,\\
L^\mathrm{ad}(X^\mathrm{ad},s,U^\mathrm{ad}) \\
= \mathbb{E}\langle P^\mathrm{ad}(0,s(0)) (X_0^\mathrm{ad} - \bar{X}_0^\mathrm{ad}), X_0^\mathrm{ad} - \bar{X}_0^\mathrm{ad} \rangle \\
+ \mathbb{E}\langle \bar{P}^\mathrm{ad}(0,s(0))  \bar{X}_0^\mathrm{ad}, \bar{X}_0^\mathrm{ad}  \rangle  + \mathbb{E}[\delta^\mathrm{ad}(0,s(0))],
\end{array}
\right.
\end{align*}
Thus, it is obtained that $P^\mathrm{rs}$ and $P^\mathrm{ad}$ solve the same equation
\begin{align*}
\begin{array}{l}
\dot{P} + Q +  P B_1+B_1^*P \\
- P  \Big[ \sum_{i \in \mathcal{I}_+} B_{2i} R_i^{-1} B_{2i}^* - 2 \lambda S_0S_0^*  \Big]  P= 0,
\end{array}
\end{align*}
and $\bar{P}^\mathrm{rs}$ and $\bar{P}^\mathrm{ad}$ solve the same equation
\begin{align*}
\begin{array}{l}
\dot{\bar{P}} + Q +\bar{Q} + \bar{P} (B_1 + \bar{B}_1) + (B_1 + \bar{B}_1)^*\bar{P} \\
-\bar{P}\Big[\sum_{i \in \mathcal{I}_+} (B_{2i}+\bar{B}_{2i})(R_i+\bar{R}_i)^{-1}(B_{2i}+\bar{B}_{2i})^*\Big] \bar{P}  = 0.
\end{array}
\end{align*}
{\color{blue} In contrast to the mean-field free cases above, mean-field-type robust  design provides an interesting feature: the expectation of the states of the two problems coincide:
 $\bar{X}^\mathrm{rs} = \bar{X}^\mathrm{ad}$, and solve the following ODE:

\begin{align*}
\begin{array}{l}
d\bar{X} = [ (B_1 + \bar{B}_1) \\
-  \sum_{i \in \mathcal{I}_+} (B_{2i}+\bar{B}_{2i}) (R_i+\bar{R}_i)^{-1} (B_{2i}+\bar{B}_{2i})^* \bar{P} ] \bar{X}dt,\\
\bar{X}(0)=\mathbb{E}[X^\mathrm{rs}_0]=\mathbb{E}[X^\mathrm{ad}_0].
\end{array}
\end{align*}
This means that when facing a risk-sensitive problem, the "true" team decision-maker can design a robust control law of mean-field type  such that $V_j^{\mathrm{ad}} 
= \sqrt{2\lambda}[-\tilde{R}_j]^{-\frac{1}{2}}S_0^*P^\mathrm{ad}(X^\mathrm{ad}-\bar{X}^\mathrm{ad}), \ j \in \mathcal{I}_{-}.$ Moreover this robust control law of mean-field type
 preserves the mean trajectory of the states because  $\mathbb{E}V_j^{\mathrm{ad}} =0.$

}

\subsection{Mean-Field-Type Game Problem}
We now discuss the relationship between robustness and risk-sensitivity in the context of games. 
Consider the risk-sensitive mean-field-type game problem given by
\begin{align*}
\begin{array}{l}
\inf_{U_i \in \mathcal{U}_i} \frac{1}{\lambda_i}\log \left[ \mathbb{E}e^{\lambda_i L_i(X,s,U)} \right],\\
\text{s.t.}~~ dX = [ B_1(X - \bar{X}) + (B_1 + \bar{B}_1)\bar{X} \\
+ \sum_{j=1}^{n} B_{2j} (U_j - \bar{U}_j) + \sum_{j=1}^{n} (B_{2j}+\bar{B}_{2j}) \bar{U}_j ]dt
+ S_0 dB,
\end{array}
\end{align*}
where $\lambda_i>0, \ L_i(X,s,U)$ is the cost presented in \eqref{eq:main_cost}. Hence, by Theorem \ref{thm:2} it follows that the solution is given by
\begin{align*}
\begin{array}{l}
U_i  =-  R_i^{-1}B_{2i}^*P^{\mathrm{rs}}_i(X - \bar{X})\\
- (R_i+\bar{R}_i)^{-1} (B_{2i}+\bar{B}_{2i})^* \bar{P}^{\mathrm{rs}}_i\bar{X},
\end{array}
\end{align*}

We design the following robust mean-field-type game problem in which each "true" decision-maker $i$ of the game above can also treat the risk as a  
{\it  fictitious adversary}  $ik$ who picks a 
robust control strategy $V_{ik}.$ 
\begin{align*}
\begin{array}{l}
\inf_{U_i \in \mathcal{U}_i} \Big[ \sup_{V_{ik}} L_i^\mathrm{ad}(X,s,U, V) \Big] \\
\text{s.t.}~~ dX = [ B_1(X - \bar{X}) + (B_1 + \bar{B}_1)\bar{X} \\
+ \sum_{i=1}^{n} B_{2i} (U_i - \bar{U}_i) + \sum_{i=1}^{n} (B_{2i}+\bar{B}_{2i}) \bar{U}_i ]dt\\
+ S_0 dB + \textcolor{blue}{\sum_{i=1}^n  \sqrt{2\lambda_i}S_0[-\tilde{R}_{ik}]^{\frac{1}{2}}(V_{ik}-\bar{V}_{ik})dt},\\
-\tilde{R}_{ik}\succ 0,
\end{array}
\end{align*}
where $L_i(X,s,U)$ is the cost presented in (\ref{eq:main_cost}), the robust cost functional 
\begin{align*}
\begin{array}{l}L_i^\mathrm{ad}(X,s,U, V) =\\
L_i(X,s,U)
+ \int_0^T \sum_{i=1}^n \langle \tilde{R}_{ik}(V_{ik} - \bar{V}_{ik}) ,V_{ik} - \bar{V}_{ik} \rangle dt\\ 
+ \int_0^T \sum_{i=1}^n \langle (\tilde{R}_{ik}+\tilde{\bar{R}}_{ik}) \bar{V}_{ik} , \bar{V}_{ik} \rangle dt,\\
-(\tilde{R}_{ik}+\tilde{\bar{R}}_{ik})\succ 0,
\end{array}
\end{align*}

Hence, by Theorem \ref{thm:2} it follows that the solution is given by:
\begin{align*}
\left\{
\begin{array}{l}
U^\mathrm{ad}_i  =-  R_i^{-1}B_{2i}^*P^\mathrm{ad}_i(X^\mathrm{ad} - \bar{X}^\mathrm{ad})\\
- (R_i+\bar{R}_i)^{-1} (B_{2i}+\bar{B}_{2i})^* \bar{P}^\mathrm{ad}_i\bar{X}^\mathrm{ad},\\
V^\mathrm{ad}_{ik}  = -\tilde{R}_{ik}^{-1}B_{2ik}^*P^\mathrm{ad}_i(X^\mathrm{ad} - \bar{X}^\mathrm{ad}),\\
\bar{V}^\mathrm{ad}_{ik}=0, 
\end{array}
\right.
\end{align*}
The coefficients $P^{\mathrm{rs}}_i$ and $P_i^\mathrm{ad}$ solve the corresponding differential equation presented next:
\begin{equation*}
\left\{
\begin{array}{l}
\dot{P}_i + Q_i +P_i  B_1^*+ B_1P_i\\
- P_i { [ B_{2i}  R_i^{-1} B_{2i}^* -2\lambda_i S_0S_0^*  ] }P_i \\
-  \sum_{j \ne i} P_j B_{2j} R_j^{-1*}B_{2j}^*P_i\\
-  \sum_{j \ne i} P_i B_{2j} R_j^{-1*}B_{2j}^*P_j = 0,\\
P_i(T,s) = Q_i(T,s),\\
\dot{\bar{P}}_i + Q_i+\bar{Q}_i + \bar{P}_i(B_1 + \bar{B}_1)^*+ (B_1 + \bar{B}_1)\bar{P}_i \\
- \bar{P}_i{\color{black}[ (B_{2i}+\bar{B}_{2i})(R_i+\bar{R}_i)^{-1}(B_{2i}+\bar{B}_{2i})^* ] }\bar{P}_i \\
-  \sum_{j \neq i} \bar{P}_j (B_{2j}+\bar{B}_{2j}) (R_j+\bar{R}_j)^{-1*} (B_{2j}+\bar{B}_{2j})^* \bar{P}_i\\
-  \sum_{j \neq i} \bar{P}_i (B_{2j}+\bar{B}_{2j}) (R_j+\bar{R}_j)^{-1*} (B_{2j}+\bar{B}_{2j})^* \bar{P}_j= 0,\\
\bar{P}_i(T,s) = Q_i(T,s)+\bar{Q}_i(T,s).
\end{array}
\right.
\end{equation*}
{\color{blue} 
Moreover, the expected matrix-value of the states $\bar{X}^{\mathrm{rs}}$ and $\bar{X}^\mathrm{ad}$ coincides. 
This result is very interesting because it says that for this class of risk-sensitive mean-field-type game problem, one can design
 a risk-neutral robust mean-field-type control strategy that preserves the mean trajectory in addition of sharing similar equilibrium control law structures.  
 The optimal control laws of the fictitious adversaries $V^\mathrm{ad}_{ik}$ satisfies $\mathbb{E}V^\mathrm{ad}_{ik}= \bar{V}^\mathrm{ad}_{ik}=0.$ 
 By choosing the following specific design 
 \begin{equation}
 \begin{array}{ll}
 S_0= \sqrt{\mu} D,\\
 -\tilde{R}_i=\gamma^2_i \mathbb{I}_{d\times d} \succ 0, \\
 \gamma_i\sqrt{2\lambda_i \mu } =1,\ 
 \end{array}  
 \end{equation} 
 the robust mean-field-type game yields }
  \begin{equation}
 \begin{array}{ll}
 L_i^\mathrm{ad}(X,s,U, V) =\\
 L_i(X,s,U)
-\gamma^2_i \int_0^T \langle (V_{ik} - \bar{V}_{ik}) ,V_{ik} - \bar{V}_{ik} \rangle dt,\\
dX = [ B_1(X - \bar{X}) + (B_1 + \bar{B}_1)\bar{X} \\
+ \sum_{i=1}^{n} B_{2i} (U_i - \bar{U}_i) + \sum_{i=1}^{n} (B_{2i}+\bar{B}_{2i}) \bar{U}_i ]dt\\
+ \sqrt{\mu} D dB + \textcolor{blue}{\sum_{i=1}^n  D(V_{ik}-\bar{V}_{ik})dt},\\
 \end{array}  
 \end{equation} 
 Hence we retrieve a modified version of the robust control design of risk-sensitive mean-field games proposed in \cite{TAC2014}. In contrast to \cite{TAC2014} where the mean trajectories  $\bar{X}^{\mathrm{rs}}$ and $\bar{X}^\mathrm{ad}$ can be different, here, thanks to the mean-field term $(V_{ik}-\bar{V}_{ik})$ the resulting mean of mean-field terms  of the two problems coincide. This new result is due to variance reduction design $\langle (V_{ik} - \bar{V}_{ik}) ,V_{ik} - \bar{V}_{ik} \rangle.$  The fictitious adversary control law can be seen as a noise (zero mean):  $\mathbb{E}V^\mathrm{ad}_{ik}=0$ by setting $${\color{blue} d\tilde{B}=  \sqrt{\mu} D dB +\sum_{i=1}^n  D(V_{ik}-\bar{V}_{ik})dt},$$
 Table \ref{tab:resultssum} below summarizes some relationships between the two problems.

{\tiny
 \begin{table}[htb!]
  	\hspace{-2cm}
  		
  		\begin{tabular}{|p{3cm}|p{5cm}||p{7cm}|}
  			\hline
  			& \multicolumn{2}{c}{\textbf{Reference}}\\
			
  			\hline
  			\textbf{RS without mean-field} &  classical RN robust control design  & this work  \\
  			\hline
			1 RS decision-maker &   -RN mean-field free robust control design            & $\checkmark$ RN mean-field type  robust control design \\
                                &  - different mean state  than the original RS problem  & $\checkmark$ same mean state as the original RS problem  \\
                                &   -designed a fictitious RN adversary                  & $\checkmark$ designed a fictitious RN adversary  \\
                                & with non-zero mean contribution to the state          & with zero mean contribution to the state  \\     
                                \hline \hline         
		    1 RS  team (full cooperation)  &   -RN mean-field free robust control design            & $\checkmark$ RN mean-field type  control design \\
                                           &  - different mean state  than the original RS problem  & $\checkmark$ same mean state as the original RS problem  \\
                                           &   -designed a fictitious RN adversary                  & $\checkmark$ designed a fictitious RN adversary  \\
                                           & with non-zero mean contribution to the state          &  with zero mean contribution to the state  \\     
                                \hline \hline         
			 n- RS non-cooperative DMs &   -RN mean-field free robust control design            & $\checkmark$ RN mean-field type  control design \\
                                       &   -different mean state  than the original RS problem  & $\checkmark$ same mean state as the original RS problem  \\
                                       &   -designed a fictitious RN adversary                  & $\checkmark$ designed a fictitious RN adversary  \\
                                       & with non-zero mean contribution to the state          &  with zero mean contribution to the state  \\     
                                \hline \hline
                                \end{tabular}
  	\caption{Contribution to Mean-Field-Free Games. RS stands for risk-sensitive, RN: risk neutral, DM: decision-maker.}
  	\label{tab:resultssum}
  \end{table}
  
  \begin{table}[htb!]
  	\hspace{-2cm}
  		
  		\begin{tabular}{|l|l|l|}
  			\hline
  			& \multicolumn{2}{c}{\textbf{Reference}}\\
			\hline
  			\textbf{RS mean-field-type dependent} & classical  & this work  \\
  			\hline
		    1 RS  decision-maker      &      & $\checkmark$ RN mean-field type  robust control design \\
                                      &      & $\checkmark$ same mean state as the original RS problem  \\
                                      &      & $\checkmark$ designed a fictitious RN adversary  \\
                                      &      &  with zero mean contribution to the state  \\     
                               \hline \hline          
			1 RS  team (full cooperation) &     & $\checkmark$ RN mean-field type  robust control design \\
                                          &     & $\checkmark$ same mean state as the original RS problem  \\
                                          &     & $\checkmark$ designed 1  fictitious RN adversary  \\
                                          &     &  with zero mean contribution to the state  \\     
                               \hline \hline
			 n- RS non-cooperative DMs &     & $\checkmark$ RN mean-field type robust  control design \\
                                       &     & $\checkmark$ same mean state as the original RS problem  \\
                                       &     & $\checkmark$ designed  n-fictitious RN adversaries  \\
                                       &     &  with zero mean contribution to the state  \\     
                               \hline \hline

						%
  		\end{tabular}
  	\caption{Contributions  to mean-field-dependent games. RS stands for risk-sensitive, RN: risk neutral, DM: decision-maker.}
  	\label{tab:resultssum}
  \end{table}
  }

\section{Proof of the Main Results}

\subsection{Proof of the risk-neutral case (Theorem \ref{thm:1})}
First, inspired from the form of the cost function consider the following guess functional:
\begin{align*}
F_i(t,s,X) = \langle P_i (X - \bar{X}), X - \bar{X}  \rangle + \langle \bar{P}_i  \bar{X}, \bar{X}  \rangle  + \delta_i.
\end{align*}

Taking the expectation of It\^o's formula applied to $F_i$ yields
\begin{equation*}
\begin{array}{l}
\mathbb{E}[dF_i] = \mathbb{E}\Big[ \langle \dot{P}_i (X - \bar{X}), X - \bar{X}  \rangle + \langle \dot{\bar{P}}_i  \bar{X}, \bar{X}  \rangle + \dot{\delta}_i\\
+ \langle (P^*_i+P_i) (X - \bar{X}),\\
  B_1(X(t) - \bar{X}(t))+ \sum_{j=1}^{n} B_{2j} (U_j - \bar{U}_j) \rangle\\
+ \langle (\bar{P}^*_i+\bar{P}_i) \bar{X} ,  (B_1 + \bar{B}_1)\bar{X}(t) + \sum_{j=1}^{n} (B_{2j}+\bar{B}_{2j}) \bar{U}_j \rangle\\
+ \frac{1}{2} \langle (P^*_i+P_i)S_0,S_0\rangle + \frac{1}{2} \int_{\Theta} \langle (P^*_i+P_i)M_0,M_0 \nu(d\theta) \rangle \Big] dt\\
+ \langle [\sum_{s'\neq s}({P}_i(t,s')-{P}_i(t,s))\tilde{q}_{ss'}] (X - \bar{X}), X - \bar{X}  \rangle dt \\
 + \langle  [\sum_{s'\neq s}({\bar{P}}_i(t,s')-\bar{P}_i(t,s))\tilde{q}_{ss'}] \bar{X}, \bar{X}  \rangle dt\\
 +[\sum_{s'\neq s}({\delta}_i(t,s')-\delta_i(t,s))\tilde{q}_{ss'}] dt
\end{array}
\end{equation*}

Now we compute the expectation of the difference $L_i - F_i(0)$ yields 
\begin{align}
\label{eq:step}
\begin{array}{l}
\mathbb{E}[(L_i - F_i(0))] \\
= \mathbb{E} \langle ( Q_i(T,s(T)) - P_i(T,s(T)) ) (X(T) - \bar{X}(T)),\\
 X(T) - \bar{X}(T)  \rangle \\
+ \mathbb{E} \langle (Q_i(T,s(T))+\bar{Q}_i(T,s(T))- \bar{P}_i(T,s(T))) \bar{X}(T), \bar{X}(T)  \rangle \\
+ \mathbb{E}[0 - \delta_i(T,s(T))] \\
+ \mathbb{E}\int_{0}^{T} \langle [ \dot{P}_i + Q_i + B_1^* (P^*_i+P_i) \\
+ \sum_{s'\neq s}({P}_i(t,s')-{P}_i(t,s))\tilde{q}_{ss'} ] (X - \bar{X}), X - \bar{X}  \rangle dt\\
+ \mathbb{E}\int_{0}^{T} \langle [\dot{\bar{P}}_i + Q_i+\bar{Q}_i + (B_1 + \bar{B}_1)^*(\bar{P}^*_i+\bar{P}_i) \\
+ \sum_{s'\neq s}({\bar{P}}_i(t,s')-\bar{P}_i(t,s))\tilde{q}_{ss'}]  \bar{X}, \bar{X}  \rangle dt\\
+ \mathbb{E}\int_{0}^{T} \dot{\delta}_i dt + \frac{1}{2} \mathbb{E}\int_{0}^{T} \langle (P^*_i+P_i)S_0,S_0\rangle dt \\
+ \frac{1}{2} \mathbb{E}\int_{0}^{T}\int_{\Theta} \langle (P^*_i+P_i)M_0,M_0 \nu(d\theta) \rangle dt\\
+\mathbb{E}\int_{0}^{T} [\sum_{s'\neq s}({\delta}_i(t,s')-\delta_i(t,s))\tilde{q}_{ss'}] dt\\
+ \mathbb{E}\int_{0}^{T} \langle R_i (U_i - \bar{U}_i), U_i - \bar{U}_i  \rangle dt \\
+ \mathbb{E}\int_{0}^{T} \langle B_{2i}^* (P^*_i+P_i) (X - \bar{X}),   U_i - \bar{U}_i \rangle dt\\
+ \mathbb{E}\int_{0}^{T} \langle (R_i+\bar{R}_i) \bar{U}_i, \bar{U}_i  \rangle dt\\
+ \mathbb{E}\int_{0}^{T} \langle (B_{2i}+\bar{B}_{2i})^* (\bar{P}^*_i+\bar{P}_i) \bar{X} , \bar{U}_i \rangle dt\\
%
+ \mathbb{E}\int_{0}^{T} \langle (P^*_i+P_i) (X - \bar{X}),  \sum_{j \ne i} B_{2j} (U_j - \bar{U}_j) \rangle dt\\
+ \mathbb{E}\int_{0}^{T} \langle (\bar{P}^*_i+\bar{P}_i) \bar{X} , \sum_{j \ne i} (B_{2j}+\bar{B}_{2j}) \bar{U}_j \rangle dt
\end{array}
\end{align}
\newpage 
Now we perform the following square completion for the terms involving $U_i - \bar{U}_i$. 
First, consider the following structure:
\begin{align*}
\begin{array}{l}
|L_1[L_2(U_i - \bar{U}_i) {\color{black}+} L_3(X - \bar{X})]|^2 \\
= \langle L_1L_2(U_i - \bar{U}_i) {\color{black}+} L_1 L_3(X - \bar{X}),\\
 L_1L_2(U_i - \bar{U}_i) {\color{black}+} L_1 L_3(X - \bar{X}) \rangle \\
= \langle L_1L_2(U_i - \bar{U}_i) , L_1L_2(U_i - \bar{U}_i) \rangle \\
{\color{black}+} 2 \langle L_1 L_3(X - \bar{X}), L_1L_2(U_i - \bar{U}_i) \rangle \\
+ \langle L_1 L_3(X - \bar{X}), L_1 L_3(X - \bar{X}) \rangle.
\end{array}
\end{align*}
Then, it is concluded that
\begin{align*}
\begin{array}{l}
\langle L_1L_2(U_i - \bar{U}_i) , L_1L_2(U_i - \bar{U}_i) \rangle \\
= \langle R_i (U_i - \bar{U}_i), U_i - \bar{U}_i  \rangle,\\
2 \langle L_1 L_3(X - \bar{X}), L_1L_2(U_i - \bar{U}_i) \rangle \\
= \langle B_{2i}^* (P^*_i+P_i) (X - \bar{X}),   U_i - \bar{U}_i \rangle.
\end{array}
\end{align*}
Thus, by setting $L_1= R_i^{-\frac{1}{2}}$, $L_2= R_i$, and $L_3 = \frac{1}{2} B_{2i}^* (P^*_i+P_i)$ yields
\begin{align*}
\begin{array}{l}
L_2^*L_1^*L_1L_2 = R_i,\\
L_2^*L_1^*L_1 L_3 = B_{2i}^* (P^*_i+P_i).
\end{array}
\end{align*}
Hence, 
\begin{align}
\label{eq:completion1}
\begin{array}{l}
\mathbb{E}\int_{0}^{T} \langle R_i (U_i - \bar{U}_i), U_i - \bar{U}_i  \rangle dt \\
+ \mathbb{E}\int_{0}^{T} \langle B_{2i}^* (P^*_i+P_i) (X - \bar{X}),   U_i - \bar{U}_i \rangle dt\\
= \mathbb{E}\int_{0}^{T} |R_i^{\frac{1}{2}}[U_i - \bar{U}_i {\color{black}+} \frac{1}{2} R_i^{-1}B_{2i}^* (P^*_i+P_i)(X - \bar{X})]|^2 dt \\
- \frac{1}{4} \mathbb{E}\int_{0}^{T} \langle (P^*_i+P_i) B_{2i} R^{-\frac{1}{2}*} R^{-\frac{1}{2}} B_{2i}^* (P^*_i+P_i)(X - \bar{X}),\\
 X - \bar{X} \rangle dt.
\end{array}
\end{align}

Using the same procedure for the terms $\bar{U}_i$ with the structure $|\tilde{L}_1[\tilde{L}_2\bar{U}_i {\color{black}+} \tilde{L}_3\bar{X}]|^2$ yields
\begin{align*}
\begin{array}{l}
\langle \tilde{L}_1\tilde{L}_2(U_i - \bar{U}_i) , \tilde{L}_1\tilde{L}_2(U_i - \bar{U}_i) \rangle \\
= \langle (R_i+\bar{R}_i) \bar{U}_i, \bar{U}_i  \rangle,\\
2 \langle \tilde{L}_1 \tilde{L}_3(X - \bar{X}), \tilde{L}_1\tilde{L}_2(U_i - \bar{U}_i) \rangle \\
= \langle (B_{2i}+\bar{B}_{2i})^* (\bar{P}^*_i+\bar{P}_i) \bar{X} , \bar{U}_i \rangle.
\end{array}
\end{align*}
Thus, by setting $\tilde{L}_1= (R_i+\bar{R}_i)^{-\frac{1}{2}}$, $\tilde{L}_2= R_i+\bar{R}_i$, and $\tilde{L}_3 = \frac{1}{2} (B_{2i}+\bar{B}_{2i})^* (\bar{P}^*_i+\bar{P}_i)$ yields
\begin{align}
\label{eq:completion2}
\begin{array}{l}
\mathbb{E}\int_{0}^{T} \langle (R_i+\bar{R}_i) \bar{U}_i, \bar{U}_i  \rangle + \langle (B_{2i}+\bar{B}_{2i})^* (\bar{P}^*_i+\bar{P}_i) \bar{X} , \bar{U}_i \rangle dt\\
= \mathbb{E}\int_{0}^{T} |(R_i+\bar{R}_i)^{\frac{1}{2}}[\bar{U}_i \\
{\color{black}+} \frac{1}{2} (R_i+\bar{R}_i)^{-1} (B_{2i}+\bar{B}_{2i})^* (\bar{P}^*_i+\bar{P}_i)\bar{X}]|^2 dt\\
- \frac{1}{4} \mathbb{E}\int_{0}^{T} \langle (\bar{P}^*_i+\bar{P}_i)(B_{2i}+\bar{B}_{2i})(R_i+\bar{R}_i)^{-\frac{1}{2}*}\\
(R_i+\bar{R}_i)^{-\frac{1}{2}}(B_{2i}+\bar{B}_{2i})^* (\bar{P}^*_i+\bar{P}_i)
\bar{X}, \bar{X} \rangle dt.
\end{array}
\end{align}
%
%

From \eqref{eq:step}, and replacing the square completions in \eqref{eq:completion1}-\eqref{eq:completion2}, we obtain that

\begin{align*}
\begin{array}{l}
\mathbb{E}[(L_i - F_i(0))] \\
= \mathbb{E} \langle ( Q_i(T,s(T)) - P_i(T,s(T)) ) (X(T) - \bar{X}(T)),\\
 X(T) - \bar{X}(T)  \rangle \\
+ \mathbb{E} \langle (Q_i(T,s(T))+\bar{Q}_i(T,s(T))- \bar{P}_i(T,s(T))) \bar{X}(T), \bar{X}(T)  \rangle \\
+ \mathbb{E}[0 - \delta_i(T,s(T))] \\
+ \mathbb{E}\int_{0}^{T} \langle [ \dot{P}_i + Q_i + B_1^* (P^*_i+P_i) \\
+ \sum_{s'\neq s}({P}_i(t,s')-{P}_i(t,s))\tilde{q}_{ss'} ] (X - \bar{X}), X - \bar{X}  \rangle dt\\
+ \mathbb{E}\int_{0}^{T} \langle [\dot{\bar{P}}_i + Q_i+\bar{Q}_i + (B_1 + \bar{B}_1)^*(\bar{P}^*_i+\bar{P}_i) \\
+ \sum_{s'\neq s}({\bar{P}}_i(t,s')-\bar{P}_i(t,s))\tilde{q}_{ss'}]  \bar{X}, \bar{X}  \rangle dt\\
+ \mathbb{E}\int_{0}^{T} \dot{\delta}_i dt + \frac{1}{2} \mathbb{E}\int_{0}^{T} \langle (P^*_i+P_i)S_0,S_0\rangle dt \\
+ \frac{1}{2} \mathbb{E}\int_{0}^{T}\int_{\Theta} \langle (P^*_i+P_i)M_0,M_0 \nu(d\theta) \rangle dt\\
+\mathbb{E}\int_{0}^{T} [\sum_{s'\neq s}({\delta}_i(t,s')-\delta_i(t,s))\tilde{q}_{ss'}] dt\\
+ \mathbb{E}\int_{0}^{T} |R_i^{\frac{1}{2}}[U_i - \bar{U}_i {\color{black}+} \frac{1}{2} R_i^{-1}B_{2i}^* (P^*_i+P_i)(X - \bar{X})]|^2 dt \\
- \frac{1}{4} \mathbb{E}\int_{0}^{T} \langle (P^*_i+P_i) B_{2i} R^{-\frac{1}{2}*} R^{-\frac{1}{2}} B_{2i}^* (P^*_i+P_i)(X - \bar{X}),\\
 X - \bar{X} \rangle dt\\
 + \mathbb{E}\int_{0}^{T} |(R_i+\bar{R}_i)^{\frac{1}{2}}[\bar{U}_i \\
 {\color{black}+} \frac{1}{2} (R_i+\bar{R}_i)^{-1} (B_{2i}+\bar{B}_{2i})^* (\bar{P}^*_i+\bar{P}_i)\bar{X}]|^2 dt\\
\end{array}
\end{align*}
\begin{align*}
\begin{array}{l} 
-  \frac{1}{4}\mathbb{E}\int_{0}^{T} \langle (\bar{P}^*_i+\bar{P}_i)(B_{2i}+\bar{B}_{2i})(R_i+\bar{R}_i)^{-\frac{1}{2}*}\\
(R_i+\bar{R}_i)^{-\frac{1}{2}}(B_{2i}+\bar{B}_{2i})^* (\bar{P}^*_i+\bar{P}_i)
\bar{X}, \bar{X} \rangle dt\\
%
+ \mathbb{E}\int_{0}^{T} \frac{1}{2} \langle (P^*_i+P_i) (X - \bar{X}),\\  
\sum_{j \ne i} B_{2j} R_j^{-1}B_{2j}^* (P^*_j+P_j)(X - \bar{X}) \rangle dt\\
+ \mathbb{E}\int_{0}^{T} \frac{1}{2} \langle (\bar{P}^*_i+\bar{P}_i) \bar{X} ,\\ 
\sum_{j \ne i} (B_{2j}+\bar{B}_{2j}) (R_j+\bar{R}_j)^{-1} (B_{2j}+\bar{B}_{2j})^* (\bar{P}^*_j+\bar{P}_j)\bar{X} \rangle dt
\end{array}
\end{align*}

By grouping the terms one arrives at the following equality:
\begin{align*}
\begin{array}{l}
\mathbb{E}[(L_i - F_i(0))] \\
= \mathbb{E} \langle ( Q_i(T,s(T)) - P_i(T,s(T)) ) (X(T) - \bar{X}(T)),\\
 X(T) - \bar{X}(T)  \rangle \\
+ \mathbb{E} \langle (Q_i(T,s(T))+\bar{Q}_i(T,s(T))- \bar{P}_i(T,s(T))) \bar{X}(T), \bar{X}(T)  \rangle \\
+ \mathbb{E}[0 - \delta_i(T,s(T))] \\
+ \mathbb{E}\int_{0}^{T} \Big\langle \Big[ \dot{P}_i + Q_i + B_1^* (P^*_i+P_i) \\
+ \sum_{s'\neq s}({P}_i(t,s')-{P}_i(t,s))\tilde{q}_{ss'} \\
- \frac{1}{4} (P^*_i+P_i) B_{2i} \textcolor{black}{R_i}^{-\frac{1}{2}*} \textcolor{black}{R_i}^{-\frac{1}{2}} B_{2i}^* (P^*_i+P_i) \\
+ \frac{1}{2} \sum_{j \ne i} (P^*_j+P_j)B_{2j} R_j^{-1*}B_{2j}^*(P^*_i+P_i)\Big] (X - \bar{X}),\\
 X - \bar{X}  \Big\rangle dt\\
 + \mathbb{E}\int_{0}^{T} \Big\langle \Big[ \dot{\bar{P}}_i + Q_i+\bar{Q}_i + (B_1 + \bar{B}_1)^*(\bar{P}^*_i+\bar{P}_i) \\
 + \sum_{s'\neq s}({\bar{P}}_i(t,s')-\bar{P}_i(t,s))\tilde{q}_{ss'} \\
 - \frac{1}{4} (\bar{P}^*_i+\bar{P}_i)(B_{2i}+\bar{B}_{2i})(R_i+\bar{R}_i)^{-\frac{1}{2}*}\\
 (R_i+\bar{R}_i)^{-\frac{1}{2}}(B_{2i}+\bar{B}_{2i})^* (\bar{P}^*_i+\bar{P}_i) \\
 + \frac{1}{2} \sum_{j \ne i} (\bar{P}^*_j+\bar{P}_j)^* (B_{2j}+\bar{B}_{2j}) (R_j+\bar{R}_j)^{-1*}\\
 (B_{2j}+\bar{B}_{2j})^* (\bar{P}^*_i+\bar{P}_i) \Big]  \bar{X}, \bar{X} \Big\rangle dt\\
\end{array}
\end{align*}
\begin{align*}
\begin{array}{l}
+ \mathbb{E}\int_{0}^{T} \Big[ \dot{\delta}_i + \frac{1}{2} \langle (P^*_i+P_i)S_0,S_0\rangle + \frac{1}{2} \int_{\Theta} \langle (P^*_i+P_i)M_0,\\
M_0 \nu(d\theta) \rangle + \sum_{s'\neq s}({\delta}_i(t,s')-\delta_i(t,s))\tilde{q}_{ss'}\Big] dt\\
+ \mathbb{E}\int_{0}^{T} |R_i^{\frac{1}{2}}[U_i - \bar{U}_i {\color{black}+} \frac{1}{2} R_i^{-1}B_{2i}^* (P^*_i+P_i)(X - \bar{X})]|^2 dt \\
+ \mathbb{E}\int_{0}^{T} |(R_i+\bar{R}_i)^{\frac{1}{2}}[\bar{U}_i \\
{\color{black}+} \frac{1}{2} (R_i+\bar{R}_i)^{-1} (B_{2i}+\bar{B}_{2i})^* (\bar{P}^*_i+\bar{P}_i)\bar{X}]|^2 dt\\
\end{array}
\end{align*}
Finally, we are making the process identification minimizing terms obtaining the announced result. \eod

\subsection{Proof of the Risk-Sensitive Statement}
We now prove Theorem \ref{thm:2}. The martingale term in the It\^o's formula is $2\lambda_i\int_0^T \langle S_0^*P_i (X-\bar{X}), dB\rangle.$ By adding and removing the term
$$2\lambda_i^2\int_0^T \langle  P_iS_0S_0^*P_i (X-\bar{X}), (X-\bar{X})\rangle dt,$$ to the  It\^o's formula we re-organize the  exponential quadratic terms.  Using that relation 
\begin{align*}
\begin{array}{l}
\mathbb{E} \Big[\exp\Big\{ 2\lambda_i\int_0^T \langle S_0^*P_i (X-\bar{X}), dB\rangle\\
- 2\lambda_i^2\int_0^T \langle  P_iS_0S_0^*P_i (X-\bar{X}), (X-\bar{X})\rangle dt\Big\}  \Big] =1,
\end{array}
\end{align*}
and matching the term in $P_i$ we arrive at the announced result. This completes the proof.

\section{Concluding Remarks}

In this article we studied risk-sensitive, robust and risk-neutral mean-field-type games in a {\it matrix-valued} linear-quadratic setting. Non-cooperation, full cooperation, adversarial  games are analyzed and their solutions are provided in a semi-explicit way.  In addition, important relationships between these problems are established. For a class of linear-quadratic risk-sensitive mean-field-type game problems driven by (matrix-valued) Brownian motions, we have designed a {\it novel}   linear-quadratic risk-neutral robust mean-field-type game problems that preserves the same mean  trajectory. The coefficients considered in this paper are random in the sense that they are parametrized with regime switching. This opens the question of extension to other type of correlated noises and other random coefficients. We leave this question for future investigation.

\bibliographystyle{unsrt}



\newpage 

\noindent {\bf Julian Barreiro-Gomez} received his B.S. degree (cum laude) in Electronics Engineering from Universidad Santo Tomas (USTA), Bogota, Colombia, in 2011. He received the MSc. degree in Electrical Engineering and the Ph.D. degree in Engineering from Universidad de Los Andes (UAndes), Bogota, Colombia, in 2013 and 2017, respectively. He received the Ph.D. degree (cum laude) in Automatic, Robotics and Computer Vision from the Technical University of Catalonia (UPC), Barcelona, Spain, in 2017; the best Ph.D. thesis in control engineering 2017 award from the Spanish National Committee of Automatic Control (CEA) and Springer; and the EECI Ph.D. Award from the European Embedded Control Institute in recognition to the best Ph.D. thesis in Europe in the field of Control for Complex and Heterogeneous Systems 2017. He received the ISA Transactions Best Paper Award 2018 in Recognition to the best paper published in the previous year. He is currently a Post-Doctoral Associate in the Learning \& Game Theory Laboratory at the New York University in Abu Dhabi (NYUAD), United Arab Emirates. His main research interests are Mean-field-type Games, Constrained Evolutionary Game Dynamics, Distributed Optimization, and Distributed Predictive Control.\\[2cm]

\noindent {\bf Tyrone E. Duncan}  received the B.E.E. degree from Rensselaer Polytechnic Institute, Troy, NY, in 1963 and the M.S. and Ph.D. degrees from Stanford University, Stanford, CA, in 1964 and 1967, respectively. He has held regular positions with the University of Michigan, Ann Arbor (1967-1971), the State University of New York, Stony Brook (1971-1974), and the University of Kansas, Lawrence (1974- present), where he is Professor of Mathematics. He has held visiting positions with the University of California, Berkeley (1969-1970), the University of Bonn, Germany (1978-1979), and Harvard University, Cambridge, MA (1979-1980), and shorter visiting positions at numerous other institutions. Dr. Duncan is a member of the editorial boards of Communications on Stochastic Analysis, and Risk and Decision Analysis and was on the editorial board of SIAM Journal on Control and Optimization (1994-2007) as an Associate Editor and a Corresponding Editor. He is a member of AMS, MAA, and SIAM.\\[2cm]

\noindent {\bf Hamidou Tembine} received the M.S. degree in applied mathematics from Ecole Polytechnique, Palaiseau, France, in 2006 and the Ph.D. degree in computer science from the University of Avignon, Avignon, France, in 2009. He holds more than 150 scientific publications including magazines, letters, journals, and conferences. He is an author of the book on Distributed Strategic Learning for Engineers (CRC Press, Taylor \& Francis 2012), and coauthor of the book Game Theory and Learning in Wireless Networks (Elsevier Academic Press). 

\end{document}